\documentclass[12pt]{amsart}
\usepackage[pdftex]{graphicx}
\pdfoutput=1
\makeatletter

\newcommand{\Rmnum}[1]{\expandafter\@slowromancap\romannumeral #1@}
\makeatother
\usepackage{amsfonts}
\usepackage{txfonts}
\usepackage{paralist}
\usepackage{color}
\usepackage{ulem}

\usepackage{enumitem}

\usepackage{tabu}                     
\usepackage{multirow}                 
\usepackage{multicol}                 
\usepackage{multirow}                
\usepackage{float}                    
\usepackage{makecell}                 
\usepackage{booktabs}                 

\usepackage{url}
\usepackage[colorlinks,  
linkcolor=blue,
anchorcolor=blue,
citecolor=blue, 
filecolor=blue,  
urlcolor=blue,]{hyperref}   
\usepackage[doipre={doi:~}]{uri}
\usepackage[capitalise]{cleveref}
\crefname{equation}{}{}
\Crefname{equation}{Equation}{Equations}
\crefname{figure}{Fig.}{Figs.}
\crefname{proposition}{Proposition}{Propositions}
\crefname{assumption}{Assumption}{Assumptions}
\usepackage{upgreek}
\usepackage[round,authoryear]{natbib}
\usepackage{bibentry}
\nobibliography*

\usepackage{graphicx}
\usepackage{float}
\usepackage{amssymb,amsmath,amsthm,mathrsfs,bm}

\usepackage{mathtools}
\usepackage{altvars,mathextras}
\usepackage[algo2e,linesnumbered,ruled]{algorithm2e}
\usepackage{verbatim}
\numberwithin{equation}{section}

\newtheorem{thm}{Theorem}[section]
\newtheorem{lemma}{Lemma}[section]
\newtheorem{remark}{Remark}[section]

\newtheorem{cor}{Corollary}[section]

\usepackage{custom_tex}

	\allowdisplaybreaks
	
\begin{document}
		
		\title{Test for high-dimensional mean vectors via the weighted $L_2$-norm}
		\author{Jianghao Li, Zhenzhen Niu, Shizhe Hong and Zhidong Bai}

		\date{\today}
		\address{School of Mathematics $\&$ Statistics, Northeast Normal University, Changchun, China}
		\email{lijh458@nenu.edu.cn}
		\address{School of Mathematics and Statistics, Shandong Normal University, Jinan, China}
		\email{niuzz@sdnu.edu.cn}
		\address{School of Statistics and Management, Shanghai University of Finance and Economics}
		\email{hong.shizhe@163.sufe.edu.cn}
		\address{School of Mathematics $\&$ Statistics, Northeast Normal University, Changchun, China}
		\email{baizd@nenu.edu.cn}
		
		\begin{abstract}
			In this paper, we propose a novel approach to test the equality of high-dimensional mean vectors of several populations via the weighted $L_2$-norm. We establish the asymptotic normality of the test statistics under the null hypothesis. We also explain theoretically why our test statistics can be highly useful in weakly dense cases when the nonzero signal in mean vectors is present.
			Furthermore, we compare the proposed test with existing tests using simulation
			results, demonstrating that the weighted $L_2$-norm-based test statistic exhibits favorable properties in terms of both size and power.
			
		\end{abstract}
		
		
		
		\maketitle

		
		\section{Introduction}
		\label{intro}
		In recent years, high-dimensional data, such as the gene expression values in microarray data, have been an important field of research in statistics, and some significant progress has been made. High-dimensional data usually refer to data dimensions that greatly exceed the sample size, presenting formidable challenges to conventional statistical techniques that may fail or perform poorly in the face of such dimensionality.  A case in point is the limitations encountered with the classical Hotelling's $T^2$ test in multivariate analysis, which becomes inapplicable due to the singularity of the sample covariance matrix in high-dimensional situations. This deficiency has inspired many new approaches and theoretical advances in various settings.
		
		In this article, our primary interest is to test the equality of several high-dimensional mean vectors with possibly unknown and unequal covariance matrices among different populations. Let $\mathbf{x}_{i 1}, \mathbf{x}_{i 2}, \ldots, \mathbf{x}_{i n_{i}}$ be the $p$-dimensional observation vectors from the $i$-th population ($i \in\{1,\ldots, k\}$). We denote the $i$-th population mean vector by $\boldsymbol{\mu}_{i}$, and the $i$-th population covariance matrix by $\boldsymbol{\Sigma}_{i}$. The hypothesis testing framework we consider is described as follows:
		\begin{align}\label{Test for Multivariate Analysis of Variance}
			\mathrm{H}_{0}:  \boldsymbol{\mu}_{1}=\cdots=\boldsymbol{\mu}_{k} \quad\text { v.s. } \quad\mathrm{H}_{1}: \text { not } \mathrm{H}_{0}.
		\end{align}
		This issue is a part of many procedures of multivariate statistical analysis, garnering considerable attention in the literature.
		For the hypothesis \eqref{Test for Multivariate Analysis of Variance}, there exist some conventional methods whose asymptotic
		properties are established in the regime where the data dimension $p$ is fixed and $n$ tends to infinity. See, e.g., the classical likelihood ratio test under the normality assumption in \cite{anderson2003introduction}, namely Wilks’ Test Statistic
		\begin{align*}
			\frac{|\mathbf{E}|^{n / 2}}{|\mathbf{H}+\mathbf{E}|^{n / 2}},
		\end{align*}
		where  $\mathbf{E}=\sum_{i=1}^{k} \sum_{j=1}^{n_{i}}\left(\mathbf{x}_{i j}- \overline{\mathbf{x}}_{i}\right)\left(\mathbf{x}_{i j}-\overline{\mathbf{x}}_{i}\right)^{\mathrm{T}}$, $\mathbf{H}=\sum_{i=1}^{k} n_{i}\left(\overline{\mathbf{x}}_{i}-\overline{\mathbf{x}}\right)\left(\overline{\mathbf{x}}_{i}-\overline{\mathbf{x}}\right)^{\mathrm{T}}$, $\overline{\mathbf{x}}_{i}=\frac{1}{n_{i}} \sum_{j=1}^{n_{i}} \mathbf{x}_{i j}$, $\overline{\mathbf{x}}=   \frac{1}{n} \sum_{i=1}^{k} n_{i} \overline{\mathbf{x}}_{i}$  and  $n=\sum_{i=1}^{k} n_{i}$. Nevertheless, in high-dimensional regime, such tests are rendered ineffective when involving determinants or inverses of empirical covariance matrices. A seminal work by \cite{bai1996effect} attempts to appeal to the attention of statisticians on the impact of high-dimensions: how to modify and enhance traditional multivariate statistical procedures when there are big data dimensions.
		
		To solve the problem of ``the effect of dimensionality", some attempts have been made for the hypothesis \eqref{Test for Multivariate Analysis of Variance} in recent years. For instance, \cite{schott2007some} extended the test statistic provided by \cite{bai1996effect} to multi-sample problem. Nevertheless, there are some limitations on Schott's test, including the assumption of the normality, the restriction on the dimensionality and the sample size, and the common covariance matrices. When these assumptions are violated, direct application of this method may result in incorrect results.
		Recently, the U-statistic-based test in \cite{chen2010two} considered for testing the mean vectors of two populations was recently extended to the one-way MANOVA by \cite{hu2017on}, which is indeed a considerable improvement over the Schott's test. \cite{he2023high} adopted a Welch-Satterthwaite $\chi^{2}  $-approximation technique to approximate the test statistic's null distribution, which was inspired by \cite{zhangjt2020simple}. However, it is not an easy task to estimate the approximate critical value in high-dimensional setting due to its intricate involvement with the $k$ population covariance matrices. For more literature dealing with testing hypotheses \eqref{Test for Multivariate Analysis of Variance}, see, e.g., \cite{cai2014high, yamada2015testing, cao2019test} and among others.
		
		
		In this paper, in our attempt and appreciation of the \cite{jiang022nonparametric} approach, we propose a weighted $L_2$-norm-based test under mild conditions for the hypothesis \eqref{Test for Multivariate Analysis of Variance}. 
		Our unified framework encompasses some existing tests, such as the $L_{2}$-norm-based test (\cite{chen2010two,hu2017on}) and the supremum-type tests (\cite{cai2014two}) and so on. Notably, these tests emerge as special cases within our proposed framework. The newly developed high-dimensional test for the hypothesis \eqref{Test for Multivariate Analysis of Variance} not only avoids specifying an explicit relationship between the dimension $p$ and the total sample size $n$, but also refrains from assuming the normality of the populations. The underlying reason is motivated by the successful ideas presented in \cite{li2012two} and \cite{zhou2017note}. Furthermore, our newly proposed high-dimensional test demonstrates distinctive power, particularly in scenarios where nonzero signals are weakly dense under the local alternative hypothesis. Extensive simulations provide a sound rationale for our theoretical results.

		To make the following discussions precise, we detail some notations. The notation $\left \| \cdot \right \| $ means the spectral norm of a matrix or the Euclidean norm of a vector. The symbol $\stackrel{\text { a.s. }} {\longrightarrow}$ means the almost sure convergence. Hereafter, we denote $\stackrel{\text { p }} {\longrightarrow}$ and $\stackrel{\text { d }} {\longrightarrow}$ as probability convergence and the convergence in distribution, respectively. Given a vector  $\boldsymbol{u} \in \mathbb{R}^{m}$ , denote by  $\|\boldsymbol{u}\|_{p}=(\sum_{i=1}^{m} \vert u_{i}\vert^{p})^{1 / p}, p \geq 1$  its  $\ell_{p}$  norm. We also denote the trace of a matrix $\mathbf{A}$ by $\operatorname{tr}(\mathbf{A})$. And we also use use the notations big-$O$ and little-$o$ to express the asymptotic relationship between two quantities. Lastly, we denote  $\frac{n !}{\left(n -r\right) !}$ by $P_{n}^{r}$.

		The remainder of this paper is structured as follows: In Section~\ref{section2.1},  we formulate our test statistic for the hypothesis \eqref{Test for Multivariate Analysis of Variance} and give the asymptotic distribution of our test statistic under the null hypothesis. 
		The asymptotic power of our test is examined in Section~\ref{section2.2}, where we also provide a brief discussion of the heteroscedastic one-way MANOVA. Extensive simulation studies are presented in 
		Section ~\ref{section3}. The major results' technical proofs are described in \ref{Appendix}. 
		
		\section{Methodology and main results
		}
		\subsection{Test statistic and its asymptotic null distribution}\label{section2.1}
		In this subsection, we will present our test procedure and the asymptotic null distribution of the proposed test statistic.
		
		To start with, we note that the equation 
		\begin{align}
			\sum_{i<l}^{k}\left\| \left( \boldsymbol{\mu}_{i} - \boldsymbol{\mu}_{l} \right)\right\|^{2} \geq 0,
		\end{align}
		which equals zero if only and only $\boldsymbol{\mu}_{1}=\cdots=\boldsymbol{\mu}_{k}$, where $\left\| \left( \boldsymbol{\mu}_{i} - \boldsymbol{\mu}_{l} \right)\right\|^{2}$ can be be estimated by the two-sample
		statistic based on the $i$-th and $l$-th samples. 
		Morever, note that
		\begin{align*}
			\boldsymbol{\mu}_{i} = \boldsymbol{\mu}_{l}  \Leftrightarrow& \boldsymbol{\alpha}^{\top} ( \boldsymbol{\mu}_{i} -\boldsymbol{\mu}_{l}) = 0, \text { for any } \boldsymbol{\alpha}  \in \mathbb{R}^{p} \\
			\Leftrightarrow& ( \boldsymbol{\mu}_{i} -\boldsymbol{\mu}_{l})^{\top} \boldsymbol{\alpha} \boldsymbol{\alpha}^{\top} ( \boldsymbol{\mu}_{i} -\boldsymbol{\mu}_{l}) = 0.
		\end{align*}
		Then, the equation that $\boldsymbol{\mu}_{1}=\cdots=\boldsymbol{\mu}_{k}$ holds 
		is equivalent to testing whether $\sum_{i<l}^{k}\left\| \boldsymbol{\alpha}^{\top} \left( \boldsymbol{\mu}_{i} - \boldsymbol{\mu}_{l} \right)\right\|^{2} =0$.
		Equivalently, we can measure 
		$\boldsymbol{\mu}_{1}=\cdots=\boldsymbol{\mu}_{k}$ by the distance $\sum_{i<l}^{k}\left( \boldsymbol{\mu}_{i}- \boldsymbol{\mu}_{l} \right)^{\top}  (\mathbf{I}_{p}  +  \boldsymbol{\alpha}\boldsymbol{\alpha}^{\top})  \left( \boldsymbol{\mu}_{i}- \boldsymbol{\mu}_{l} \right)$. To adopt a more comprehensive framework, we further propose to measure the null hypothesis \eqref{Test for Multivariate Analysis of Variance} with
		\begin{align}\label{Expecation of test statistic}
			\sum_{i<l}^{k}\left( \boldsymbol{\mu}_{i}- \boldsymbol{\mu}_{l} \right)^{\top}  ( \mathbf {B}  +  \boldsymbol{\alpha}\boldsymbol{\alpha}^{\top})  \left( \boldsymbol{\mu}_{i}- \boldsymbol{\mu}_{l} \right),
		\end{align}
		where $\mathbf {B} = \operatorname{diag} \left\{\omega_{1}^{2}, \omega_{2}^{2}, \ldots, \omega_{p}^{2}\right\}, \boldsymbol{\alpha} =\left( \alpha_{1}, \ldots, \alpha_{p}\right)^{\top}$.
		\begin{remark}
			It follows that the first term in \eqref{Expecation of test statistic} is a weighted $L_{2}$-norm between means, and the second term can be understood as reducing to the dimension while retaining the significant signal or information present in the mean vectors.
		\end{remark}
		Denote $\mathbf{W} = \mathbf {B} +\boldsymbol{\alpha}\boldsymbol{\alpha}^{\top}$, a test statistic for the hypothesis \eqref{Test for Multivariate Analysis of Variance} can then be defined as
		\begin{align}
			T_{n}  
			=&\sum_{i<l}^{k}
			\left( \bar{\mathbf{x}}_{i}- \bar{\mathbf{x}}_{l} \right)^{\top} \mathbf{W}\left( \bar{\mathbf{x}}_{i}- \bar{\mathbf{x}}_{l}\right)-(k-1) \sum_{i=1}^{k} n_{i}^{-1} \operatorname{tr} \mathbf{W} \mathbf{S}_{i} \notag  \\
			=&(k-1) \sum_{i=1}^{k} \frac{1}{n_{i}\left(n_{i}-1\right)} \sum_{j_{1} \neq j_{2}}^{n_{i}} \mathbf{x}_{i j_{1}}^{\top}  \mathbf{W} \mathbf{x}_{i j_{2}}-\sum_{i<l}^{k} \frac{2}{n_{i} n_{l}} \sum_{j_{1}}^{n_{i}} \sum_{j_{2}}^{n_{l}} \mathbf{x}_{i j_{1}}^{\top}  \mathbf{W} \mathbf{x}_{l j_{2}}.
		\end{align}
		Elementary derivation shows 
		\begin{align}
			\mathbb{E} \left( T_{n} \right) =& \sum_{i<l}^{k}\left( \boldsymbol{\mu}_{i}- \boldsymbol{\mu}_{l} \right)^{\top} \mathbf{W}\left( \boldsymbol{\mu}_{i}- \boldsymbol{\mu}_{l} \right) \notag \\
			=& \sum_{i<l}^{k}\left( \boldsymbol{\mu}_{i}- \boldsymbol{\mu}_{l} \right)^{\top}  \mathbf{B}\left( \boldsymbol{\mu}_{i}- \boldsymbol{\mu}_{l} \right) + \sum_{i<l}^{k}\left( \boldsymbol{\mu}_{i}- \boldsymbol{\mu}_{l} \right)^{\top}  \boldsymbol{\alpha}\boldsymbol{\alpha}^{\top}\left( \boldsymbol{\mu}_{i}- \boldsymbol{\mu}_{l} \right)
		\end{align}
		and
		\begin{align*}
			\operatorname{Var}\left(T_n\right)= & \sum_{i=1}^k \frac{2(k-1)^2}{n_i\left(n_i-1\right)} \operatorname{tr}\left( \mathbf{W} \boldsymbol{\Sigma}_i\right)^2 +
			\sum_{i<l}^k \frac{4}{n_i n_l} \operatorname{tr}\left( \mathbf{W} \boldsymbol{\Sigma}_i \mathbf{W} \boldsymbol{\Sigma}_l\right) \\
			& +4 \sum_{i=1}^k \frac{1}{n_i}\left(\sum_{l=1}^k \boldsymbol{\mu}_l-k \boldsymbol{\mu}_i\right)^{\top} \mathbf{W} \boldsymbol{\Sigma}_i \mathbf{W} \left(\sum_{l=1}^k \boldsymbol{\mu}_l-k \boldsymbol{\mu}_i\right). 
		\end{align*}
		In what follows, to analyze the asymptotic behaviors of the considered statistics throughout the paper, we impose the following assumptions, which are commonly used by \cite{bai1996effect, chen2010two} and among others.
		
		\textbf{Assumption A:} There exist a  $p \times m$  matrix  $\boldsymbol{\Gamma}_{i}$ and a collection of $m$-dimensional random vectors  $\left\{\mathbf{z}_{ij}\right\}_{j=1}^{n_i}$ such that
		\begin{align} \label{general factor model}
			\mathbf{x}_{i j}=\boldsymbol{\mu}_i + \boldsymbol{\Gamma}_{i} \mathbf{z}_{i j}, i \in\{1,\ldots, k\}, j= 1, \ldots, n_i,
		\end{align}
		where $\boldsymbol{\Gamma}_{i} \boldsymbol{\Gamma}_{i}^{\top}=\boldsymbol{\Sigma}_{i}$ with $m \geq p$, and the random vectors $\mathbf{z}_{i j}=\left( z_{ij1}, \ldots, z_{ijm   }\right)^{\top}$ are independent and identically distributed (i.i.d.) with $\mathbb{E}\left( z_{i j \theta  }\right)=0, \mathbb{E}\left( z_{i  j \theta  }^2\right)=1$ and $\mathbb{E}\left( z_{i  j \theta }^4 \right)<\infty, \theta  \in \{1,\ldots, m \}$. Furthermore, given a positive integer  $a $ and  $\varsigma_{\theta }$  such that  $\sum_{\theta =1}^{a} \varsigma_{\theta } \leq 8$, the following condition is satisfied:
		\begin{align}\label{pseudo-independence}
			\mathbb{E}\left(z_{i j \theta _{1}}^{\zeta_{1}} \cdots z_{i j \theta _{a}}^{\zeta_ a}\right)=\mathbb{E} \left(z_{i j \theta _{1}}^{\zeta_1}\right) \cdots \mathbb{E} \left(z_{i j \theta _{a}}^{\zeta_a}\right), 
		\end{align}
		whenever  $\theta _{1}, \theta _{2}, \ldots, \theta _{t}$  are distinct indices.
		
		\textbf{Assumption B:} $\frac{n_{i}}{n} \rightarrow k_{i} \in (0,1), i=1, \ldots k$, as  $n \rightarrow \infty$ . Here  $n=\sum_{i=1}^{k} n_{i}$.
		
		\textbf{Assumption C:} 	
		\begin{align}\label{bounded but allow spike}
			\operatorname{tr}\left( \mathbf{W} \boldsymbol{\Sigma}_{i_1} \mathbf{W} \boldsymbol{\Sigma}_{i_2} \mathbf{W} \boldsymbol{\Sigma}_{i_3} \mathbf{W} \boldsymbol{\Sigma}_{i_4}\right)=
			o\left\{\operatorname{tr}^{2} \left[\left(\sum_{i=1}^{k} \mathbf{W} \boldsymbol{\Sigma}_{i}\right)^{2}\right]\right\}, 
		\end{align} 
		for $i_1, i_2, i_3, i_4 \in\{1,2, \ldots, k\} $.
		
		\textbf{Assumption D:} \begin{align}\left(\boldsymbol{\mu}_{i_{1}}-\boldsymbol{\mu}_{i_{2}}\right)^{\top} \mathbf{W} \boldsymbol{\Sigma}_{i_{1}} \mathbf{W} \left(\boldsymbol{\mu}_{i_{1}}-\boldsymbol{\mu}_{i_{3}}\right)=o\left\{n^{-1} \operatorname{tr}\left[\left(\sum_{i=1}^{k} \mathbf{W} \boldsymbol{\Sigma}_{i}\right)^{2}\right]\right\},
		\end{align}
		for $i_{1}, i_{2}, i_{3} \in\{1,2, \ldots, k\}$.
		
		\begin{remark}\label{assumption}
			Assumption A closely aligns with conditions (3.1) and (3.2) in \cite{chen2010two}. The equation \eqref{general factor model} is a factor model advocated earlier by \cite{bai1996effect} for testing the equality of two sample mean vectors and employed extensively in other literature. The equation \eqref{pseudo-independence} can be viewed as a kind of pseudo-independence condition of $\mathbf{z}_{i j}$, and the equation will be satisfied if $\mathbf{z}_{i j}$ has independent components. Assumption B is a standard condition in multi-sample asymptotic analysis, while Assumption C serves as a technical condition for establishing the asymptotic normality of the test statistics.  Assumption D can be understood as the local alternative hypothesis. 
			
			To grasp \eqref{bounded but allow spike},  we consider a simple case that $\lambda_{1} \leq \lambda_{2} \cdots \leq   \lambda_{p}$  and  $\lambda_{1}^{*} \leq \lambda_{2}^{*} \cdots \leq \lambda_{p}^{*}$  are eigenvalues of  $\mathbf{W}$  and  $\boldsymbol{\Sigma}_{1} = \boldsymbol{\Sigma}_{2}= \cdots=\boldsymbol{\Sigma}_{q}=\boldsymbol{\Sigma}$, respectively. Moreover, let $\alpha_{1} = \cdots = \alpha_{p}=\alpha$  and  $0<\omega_{1} \leq \cdots \leq \omega_{p}$. Then \eqref{bounded but allow spike} is equivalent to
			\begin{align*}
				\operatorname{tr}\left(\mathbf{W}\boldsymbol{\Sigma}\right)^4=
				o\left[\operatorname{tr}^{2}\left(\mathbf{W} \boldsymbol{\Sigma}\right)^{2}\right],
			\end{align*}
			and the characteristic equation of $\mathbf{W}$ is
			\begin{align*}
				f(\lambda)=\prod_{i=1}^p\left(\omega_i^2-\lambda\right)+\alpha^2\sum_{i=1}^p\prod_{j\neq i}\left(\omega_j^2-\lambda\right),
			\end{align*}
			which indicates that
			$\omega_1^2\leq\lambda_1\leq\omega_2^2,\dots,\omega_{p-1}^2\leq\lambda_{p-1}\leq\omega_p^2$. Thus, using the inequalities
			\begin{align*}
				&\operatorname{tr}\left(\mathbf{W}\boldsymbol{\Sigma}\right)^4\leq\sum_{i=1}^p\left(\lambda_i\lambda_i^*\right)^4\leq\omega_p^8\operatorname{tr}\boldsymbol{\Sigma}^4+\left(\lambda_p^4-\omega_p^8\right)\lambda_p^{*4},\\
				&\operatorname{tr}\left(\mathbf{W}\boldsymbol{\Sigma}\right)^2\geq\lambda_1^2\operatorname{tr}\boldsymbol{\Sigma}^2\geq\omega_1^4\operatorname{tr}\boldsymbol{\Sigma}^2,\\
				&\left|\lambda_p^4-\omega_p^8\right|\leq\|\mathbf{W}^4-\mathbf{B}^4\|\leq 6\left[\alpha^2\sum_{i=1}^p\omega_i^6+p\alpha^4\sum_{i=1}^p\omega_i^4+p^2\alpha^6\sum_{i=1}^p\omega_i^2+p^4\alpha^8\right],
			\end{align*}
			we have
			\begin{align*}
				\frac{\operatorname{tr}\left(\mathbf{W}\boldsymbol{\Sigma}\right)^4}{\operatorname{tr}^2\left(\mathbf{W}\boldsymbol{\Sigma}\right)^2}\leq\frac{\omega_p^8\operatorname{tr}\boldsymbol{\Sigma}^4}{\omega_1^8\operatorname{tr}^2\boldsymbol{\Sigma}^2}+6\frac{\lambda_p^{*4}\left[\alpha^2\sum_{i=1}^p\omega_i^6+p\alpha^4\sum_{i=1}^p\omega_i^4+p^2\alpha^6\sum_{i=1}^p\omega_i^2+p^4\alpha^8\right]}{\omega_1^8\operatorname{tr}^2\boldsymbol{\Sigma}^2}.
			\end{align*}
			Therefore, if $0 < \inf_{i}\omega_{i} \leq \sup_{i}\omega_i < \infty$, $\alpha^2=O\left(p^{-3/4}\right)$, $\operatorname{tr}\boldsymbol{\Sigma}^4=o\left(\operatorname{tr}^2\boldsymbol{\Sigma}^2\right)$ and $p^{1/2}\lambda_p^{*2}=o\left(\operatorname{tr}\boldsymbol{\Sigma}^2\right)$, then \eqref{bounded but allow spike} is valid. Notice that the restrictions on $\boldsymbol{\Sigma}$ only appear in the last two conditions. Similar requirements are used in \cite{chen2010two,chen2010tests,wang15}.
		\end{remark}
		\begin{thm}\label{CLT of main results}
			Under Assumptions A-D, as $p$, $n \to \infty$, we have
			\begin{align}
				\frac{T_{n }-  \sum_{i<l}^{k}\left( \boldsymbol{\mu}_{i}- \boldsymbol{\mu}_{l} \right)^{\top} \mathbf{W}\left( \boldsymbol{\mu}_{i}- \boldsymbol{\mu}_{l} \right)  }{\sigma_{n,k}}	\stackrel{\text { d }}{\longrightarrow} \mathrm{N} (0,1),
			\end{align}
			where
			\begin{align}
				\sigma_{n,k}^{2} =	\sum_{i=1}^k \frac{2(k-1)^2}{n_i\left(n_i-1\right)} \operatorname{tr}\left( \mathbf{W} \boldsymbol{\Sigma}_i\right)^2 +
				\sum_{i<l}^k \frac{4}{n_i n_l} \operatorname{tr}\left( \mathbf{W} \boldsymbol{\Sigma}_i \mathbf{W} \boldsymbol{\Sigma}_l\right) 
			\end{align}	
		\end{thm}
		
		\begin{remark}
			By Assumption D, we have $\operatorname{Var}\left(T_n\right) = \sigma_{n,k}^{2} (1+o(1))$.
		\end{remark}
		
		In particular, under the null hypothesis \eqref{Test for Multivariate Analysis of Variance}, we have the following corollary:
		\begin{cor}\label{Corollary1}
			Under Assumptions A-D and $ \mathrm{H}_{0}:  \boldsymbol{\mu}_{1}=\cdots=\boldsymbol{\mu}_{k}$, as $p$, $n \to \infty$, we have
			\begin{align}
				\frac{T_{n } }{\sigma_{n,k}}	\stackrel{\text { d }}{\longrightarrow} \mathrm{N} (0,1),
			\end{align}	
		\end{cor}
		In order to develop our test procedure, we aim to provide an asymptotically ratio-consistent estimator of $\sigma_{n,k}^{2}$. 
		In the spirit of \cite{li2012two}, we propose estimators of $ \operatorname{tr}\left(\mathbf{W} \boldsymbol{\Sigma}_{i}\right)^{2}$ and  $\operatorname{tr}\left(\mathbf{W} \boldsymbol{\Sigma}_{i} \mathbf{W} \boldsymbol{\Sigma}_{l}\right)$, and give the ratio-consistent estimator of $\sigma_{n,k}^{2}$ which is stated in the following lemma:
		\begin{lemma} \label{estimators of trace}
			Under Assumptions A-D, as $p$, $n \to \infty$, we have
			\begin{align}
				\frac{\hat{\sigma}_{n, k}^{2}}{\sigma_{n,k}^{2}} \stackrel{\text { p }}{\longrightarrow}	1,
			\end{align}
			where
			\begin{align}
				\hat{\sigma}_{n, k}^{2}=&
				\sum_{i<l}^k \frac{4}{n_i n_l} \operatorname{tr} \widehat{\left(\mathbf{W} \boldsymbol{\Sigma}_{i} \mathbf{W} \boldsymbol{\Sigma}_{l}\right)} + \sum_{i=1}^k \frac{2(k-1)^2}{n_i\left(n_i-1\right)} \operatorname{tr} \widehat{\left(\mathbf{W} \boldsymbol{\Sigma}_{i}\right)^{2}}, \notag\\
				\operatorname{tr} \widehat{ \left(\mathbf{W} \boldsymbol{\Sigma}_{i}\right)^{2} } 
				=&\frac{1}{P_{n_{i}}^{2}} \sum_{j_{1} \neq j_{2} } \mathbf{x}_{i j_{1} }^{\top}  \mathbf{W} \mathbf{x}_{i j_{2}} \mathbf{x}_{i j_{2} }^{\top}  \mathbf{W} \mathbf{x}_{i j_{1}} -\frac{2}{P_{n_{i}}^{3} }  \sum_{j_{1} \neq j_{2} \neq j_{3} } \mathbf{x}_{i j_{1}}^{\top} \mathbf{W} \mathbf{x}_{i j_{2}} \mathbf{x}_{i j_{3}}^{\top} \mathbf{W}  \mathbf{x}_{i j_{1}}    \notag\\ &+\frac{1}{P_{n_{i}}^{4} }  
				\sum_{j_{1}\neq j_{2} \neq j_{3} \neq j_{4}} \mathbf{x}_{i j_{1}}^{\top} \mathbf{W} \mathbf{x}_{i j_{2}} \mathbf{x}_{i j_{3}}^{\top} \mathbf{W} \mathbf{x}_{i j_{4}} , \\
				\operatorname{tr} \widehat{\left(\mathbf{W} \boldsymbol{\Sigma}_{i} \mathbf{W} \boldsymbol{\Sigma}_{l}\right)} 
				= &\frac{1}{n_{i} n_{l}} \sum_{j_{1}=1}^{n_{i}} \sum_{j_{2}=1}^{n_{l}} \left(\mathbf{x}_{i j_{1}}^{\top} \mathbf{W} \mathbf{x}_{ l j_{2}}\right)^{2} - \frac{1}{n_{i} n_{l}\left(n_{i}-1\right)} \sum_{j_{1}\neq j_{3}} \sum_{j_{2}}  \mathbf{x}_{l j_{2}}^{\top} \mathbf{W} \mathbf{x}_{i j_{3}}\mathbf{x}_{i j_{1}}^{\top} \mathbf{W} \mathbf{x}_{l j_{2}} \notag\\
				& -\frac{1}{n_{i} n_{l}\left(n_{l}-1\right)} \sum_{j_{2} \neq j_{4}} \sum_{j_{1}} \mathbf{x}_{i j_{1}}^{\top}\mathbf{W} \mathbf{x}_{l j_{4}} \mathbf{x}_{l j_{2}}^{\top} \mathbf{W}\mathbf{x}_{i j_{1}}  \notag \\
				& +\frac{1}{n_{i} n_{l}\left(n_{i}-1\right)\left(n_{l}-1\right)} \sum_{j_{1}\neq j_{3} }^{n_{i}} \sum_{j_{2}\neq  j_{4}}^{n_{l}} \mathbf{x}_{i j_{1}}^{\top} \mathbf{W} \mathbf{x}_{i j_{3}} \mathbf{x}_{l j_{2}}^{\top}\mathbf{W} \mathbf{x}_{l j_{4}} .
			\end{align} 
		\end{lemma}
		As highlighted by \cite{zhou2017note}, due to the translation-invariant property of U-statistics, the three distinct versions of the unbiased estimator 
		$ \operatorname{tr}\boldsymbol{\Sigma}_{i}^{2}$ that are shown in \cite{li2012two,himeno2014wstimations,hu2017on} are equivalent. Similar to \cite{zhou2017note}, by performing straightforward calculations, we can simplify the above estimators to a more computationally efficient form:
		\begin{align*}
			&\operatorname{tr} \widehat{ \left(\mathbf{W} \boldsymbol{\Sigma}_{i}\right)^{2} } \\
			=&\frac{- 1}{(n_{i} - 1) (n_{i} - 2) } \sum_{j=1 }^{n_{i}} \left( ( \mathbf{x}_{ij} -  \bar{\mathbf{x}}_{i} )^{\top} \mathbf{W} ( \mathbf{x}_{ij} -  \bar{\mathbf{x}}_{i} )\right)^{2} +\frac{1}{P_{n_{i}}^{4} } 
			(  \sum_{j=1}^{n_{i}}( \mathbf{x}_{ij} -  \bar{\mathbf{x}}_{i} )^{\top} \mathbf{W} ( \mathbf{x}_{ij} -  \bar{\mathbf{x}}_{i} )  )^{2} \\
			& + \frac{ 1}{n_{i}  (n_{i} - 3) }     \sum_{j_{1} , j_{2} }^{n_{i}} \left( ( \mathbf{x}_{ij_{1}} -  \bar{\mathbf{x}}_{i} )^{\top} \mathbf{W} ( \mathbf{x}_{ij_{2}} -  \bar{\mathbf{x}}_{i} )\right)^{2}, \\
			=& \frac{- 1}{(n_{i} - 1) (n_{i} - 2) } \sum_{j=1 }^{n_{i}} \left( ( \mathbf{x}_{ij} -  \bar{\mathbf{x}}_{i} )^{\top} \mathbf{W} ( \mathbf{x}_{ij} -  \bar{\mathbf{x}}_{i} )\right)^{2} + \frac{ (n_{i} - 1)^{2}}{n_{i}  (n_{i} - 3) }\operatorname{tr} (\mathbf{W} \mathbf{S}_{i})^{2} \\
			&+ \frac{(n_{i} - 1)}{n_{i}(n_{i} - 2) (n_{i} - 3) } \operatorname{tr}^{2} (\mathbf{W} \mathbf{S}_{i})
		\end{align*}
		and
		\begin{align*}
			\operatorname{tr} \widehat{\left(\mathbf{W} \boldsymbol{\Sigma}_{i} \mathbf{W} \boldsymbol{\Sigma}_{l}\right)} 
			= & \operatorname{tr} (\mathbf{W} \mathbf{S}_{i} \mathbf{W} \mathbf{S}_{l}  ),
		\end{align*}
		where 
		\begin{align*}
			\mathbf{S}_{i} = \frac{1}{n_{i}} \sum_{j=1}^{n_{i}} (\mathbf{x}_{ij} - \bar{\mathbf{x}}_{i})(\mathbf{x}_{ij} - \bar{\mathbf{x}}_{i})^{\top}.
		\end{align*}
		Furthermore, from Lemma \ref{estimators of trace}, we can directly derive the following theorem:
		\begin{thm}\label{Corollary of CLT of main results}
			Under Assumptions A-D, and $ \mathrm{H}_{0}:  \boldsymbol{\mu}_{1}=\cdots=\boldsymbol{\mu}_{k}$ as $p$, $n \to \infty$, we have
			\begin{align}
				\frac{T_{n } }{\hat{\sigma}_{n, k}}	\stackrel{\text { d }}{\longrightarrow} \mathrm{N} (0,1),
			\end{align}	
		\end{thm}
		According to Theorem \ref{Corollary of CLT of main results}, our proposed test rejects $\mathrm{H}_{0}$ if $T_{n } \geq \hat{\sigma}_{n, k} z_{\vartheta}$ for a given significance level $\vartheta$, where  $z_{\vartheta}$  is the upper-$\vartheta$  quantile of $N(0,1)$. 
		\subsection{Power of the proposed weighted $L_2$-norm test}\label{section2.2}
		In this subsection, we explore the power of the proposed weighted  $L_2$-norm test by providing a rough analysis of our test statistic under the local alternative hypothesis. The theorem below shows the asymptotic power function of the proposed test:
		\begin{thm}\label{Power of Corollary of CLT of main results}
			Suppose that Assumptions A-D hold. Then under $ \mathrm{H}_{1}:  \boldsymbol{\mu}_{1} \neq \cdots 
			\neq \boldsymbol{\mu}_{k} $, the power of our proposed test is given by
			\begin{align*}
				& \lim _{ n, p \rightarrow \infty} \mathrm{P}\left( T_{n } \geq \hat{\sigma}_{n, k} z_{\vartheta} \right) \\
				= & \lim _{ n, p \rightarrow \infty} \Phi\left\{-z_{\vartheta}+\frac{ \sum_{i<l}^{k}\left( \boldsymbol{\mu}_{i}- \boldsymbol{\mu}_{l} \right)^{\top} \mathbf{W}\left( \boldsymbol{\mu}_{i}- \boldsymbol{\mu}_{l} \right) }{ \hat{\sigma}_{n, k}  }\right\},
			\end{align*}
			where $ \Phi(\cdot)$ is the cumulative distribution function of the standard normal random variable.
		\end{thm}
		
		Theorem \ref{Power of Corollary of CLT of main results} provides a general result on the power of the proposed test statistic. However, analyzing the asymptotic power functions of the suggested test similar to two populations is challenging. To facilitate a basic understanding, we consider the scenario of equal population covariance matrices and introduce the following corollary:
		\begin{cor}\label{special case: Power of Corollary of CLT of main results}
			Assume $\boldsymbol{\Sigma}_{1} = \boldsymbol{\Sigma}_{2}= \cdots=\boldsymbol{\Sigma}_{k}$, and Assumptions A-D hold. Under $ \mathrm{H}_{1}:  \boldsymbol{\mu}_{1} \neq \cdots 
			\neq \boldsymbol{\mu}_{k} $ as $p$, $n \to \infty$, the power of our proposed test is given by
			\begin{align*}
				& \lim _{ n, p \rightarrow \infty} \mathrm{P}\left( T_{n } \geq \hat{\sigma}_{n, k} z_{\vartheta} \right) \\
				= & \lim _{ n, p \rightarrow \infty} \Phi\left\{-z_{\vartheta}+\frac{ \sum_{i<l}^{k}\left( \boldsymbol{\mu}_{i}- \boldsymbol{\mu}_{l} \right)^{\top} \mathbf{W}\left( \boldsymbol{\mu}_{i}- \boldsymbol{\mu}_{l} \right) }{\sqrt{2 \operatorname{tr}\left( \mathbf{W} \boldsymbol{\Sigma}_1 \right)^2  }  \sqrt{\sum_{i=1}^k \frac{(k-1)^2}{n_i\left(n_i-1\right)}  +
						\sum_{i\neq l}^k \frac{1}{n_i n_l}  } }\right\},
			\end{align*}
		\end{cor}
		Based on Corollary \ref{special case: Power of Corollary of CLT of main results}, we proceed to analyze the performance of our test in terms of the power. In alignment with the notations defined in Remark \ref{assumption}, we still use  $\lambda_{1} \leq \lambda_{2} \cdots \leq   \lambda_{p}$  and  $\lambda_{1}^{*} \leq \lambda_{2}^{*} \cdots \leq \lambda_{p}^{*}$  to represent the eigenvalues of  $\mathbf{W}$  and  $\boldsymbol{\Sigma}_{1} $ ( $= \boldsymbol{\Sigma}_{2}= \cdots=\boldsymbol{\Sigma}_{q}$), respectively. To showcase the superiority of our proposed test, we first conduct a comparative analysis with the test put forth by\cite{hu2017on}, which is considered a specific instance within our comprehensive framework. To this end, we employ asymptotic relative efficiency (ARE), a widely accepted metric for comparing the asymptotic power functions under the local alternatives. Note that the power function in \cite{hu2017on} is given by
		\begin{align*}
			\beta_{T_{hb}}
			=  \lim _{ n, p \rightarrow \infty} \Phi\left\{-z_{\vartheta}+\frac{ \sum_{i<l}^{k}\left\|\boldsymbol{\mu}_{i}-\boldsymbol{\mu}_{l}\right\|^{2}  }{\sqrt{2 \operatorname{tr}\left( \boldsymbol{\Sigma}_1 \right)^2  }  \sqrt{\sum_{i=1}^k \frac{(k-1)^2}{n_i\left(n_i-1\right)}  +
					\sum_{i\neq l}^k \frac{1}{n_i n_l}  } }\right\}.
		\end{align*}
		Then the ARE of the proposed test in comparison to the test suggested by \cite{hu2017on} is obtained as follows:
		\begin{align*}
			\operatorname{ARE}(\beta_{T_{n}}, \beta_{T_{hb}}) = \frac{ \sum_{i<l}^{k}\left( \boldsymbol{\mu}_{i}- \boldsymbol{\mu}_{l} \right)^{\top} \mathbf{W}\left( \boldsymbol{\mu}_{i}- \boldsymbol{\mu}_{l} \right) \sqrt{ \operatorname{tr}\left( \boldsymbol{\Sigma}_1 \right)^2  }             }{\sum_{i<l}^{k}\left\|\boldsymbol{\mu}_{i}-\boldsymbol{\mu}_{l}\right\|^{2}         \sqrt{ \operatorname{tr}\left( \mathbf{W} \boldsymbol{\Sigma}_1 \right)^2  }    }.
		\end{align*}
		For convenience of exposition, assume that  $\alpha_{1}=\cdots=\alpha_{p}=\alpha$  and  $\mathbf{B}=\mathbf{I}_{p}$. Under these assumptions, it can be obtained that  $\lambda_{1}=\cdots= \lambda_{p-1}=1$ and  $\lambda_{p}=1+p \alpha^{2} $. Denoting $\boldsymbol{\mu}_{i}=\left(\mu_{i1},\dots,\mu_{ip}\right)^{\top}$ for $i=1,\dots,p$, we deduce that
		\begin{align*}
			\sum_{i<l}^{k}\left(\boldsymbol{\mu}_{i}-\boldsymbol{\mu}_{l}\right)^{\top}  \mathbf{W} \left(\boldsymbol{\mu}_{i}-\boldsymbol{\mu}_{l}\right)
			=\sum_{i<l}^{k} \alpha^{2}\left(\sum_{q=1}^{p}\left(\mu_{i q}-\mu_{l q}\right)\right)^{2}+ \sum_{i<l}^{k} \sum_{q=1}^{p} \left(\mu_{i q}-\mu_{l q}\right)^{2}.
		\end{align*}
		This, together with
		\begin{align*}
			\operatorname{tr}\left( \mathbf{W} \boldsymbol{\Sigma}_1 \right)^2 \leq \sum_{i=1}^{p}\left(\lambda_{i} \lambda_{i}^{*}\right)^{2} 
			=\sum_{i=1}^{p-1} ( \lambda_{i}^{*} )^{2} +(\lambda_{p}\lambda_{p}^{*})^{2}=\operatorname{tr}\boldsymbol{\Sigma}_1^2+\left(2p\alpha^2+p^2\alpha^4\right)\lambda_{p}^{*2},
		\end{align*}
		implies that
		\begin{align*}
			\operatorname{ARE}(\beta_{T_{n}}, \beta_{T_{hb}}) \geq \frac{ 1  + \frac{\sum_{i<l}^{k} \alpha^{2}\left(\sum_{q=1}^{p}\left(\mu_{i q}-\mu_{l q}\right)\right)^{2} }{ \sum_{i<l}^{k}\left\|\boldsymbol{\mu}_{i}-\boldsymbol{\mu}_{l}\right\|^{2}               }            }{\sqrt{1 + \frac{ \left(2p\alpha^2+p^2\alpha^4\right)\lambda_{p}^{*2} }{ \operatorname{tr}\left(  \boldsymbol{\Sigma}_1 \right)^2  }  }}.
		\end{align*}
		Hence, if $\lambda_{p}^{*2}\max \left\{ p \alpha^{2},  p^{2} \alpha^{4}  \right\} = o\left( \operatorname{tr}  \boldsymbol{\Sigma}_1^2 \right) $, we arrive at the conclusion that
		\begin{align*}
			\lim _{ n, p \rightarrow \infty} \operatorname{A R E}\left(\beta_{T_{n}}, \beta_{T_{hb}}\right) \geq 1.
		\end{align*}
		This illustrates that, in an asymptotic sense, our test outperforms that proposed by \cite{hu2017on} in terms of power.

Next, we investigate the lower bound of the power for our proposed test by considering a simple case where there exists at least one pair of $i\neq j$ such that
\begin{align}\label{alternative}
\boldsymbol{\mu}_{i} - \boldsymbol{\mu}_{j} =(\underbrace{\nu, \ldots, \nu}_{p^{\delta}}, \underbrace{0, \ldots, 0}_{ p^{1-\delta}})^{\top},
\end{align} 
where $0\leq\delta\leq1$. This stands for that $\boldsymbol{\mu}_i-\boldsymbol{\mu}_j$ has $p^{\delta}$ nonzero entries of equal value. When $\alpha_{1}=\cdots=\alpha_{p}=\alpha$ and $\omega_{1} \leq   \cdots \leq \omega_{p}$, from the analysis in Remark \ref{assumption}, we obtain that
$\omega_1^2\leq\lambda_1\leq\omega_2^2,\dots,\omega_{p-1}^2\leq\lambda_{p-1}\leq\omega_p^2$. For the bound of $\lambda_p$, we have
\begin{align*}
\left|\lambda_p-\omega_p^2\right|\leq\left\|\mathbf{W}-\mathbf{B}\right\|=p\alpha^2.
\end{align*}
Thus by applying $\operatorname{tr}\left( \mathbf{W} \boldsymbol{\Sigma}_1 \right)^2\leq\lambda_p^{*2}\operatorname{tr}\mathbf{W}^2$, it's not difficult to get
\begin{align*}
&\frac{ \sum_{i<l}^{k}\left( \boldsymbol{\mu}_{i}- \boldsymbol{\mu}_{l} \right)^{\top} \mathbf{W}\left( \boldsymbol{\mu}_{i}- \boldsymbol{\mu}_{l} \right) }{\sqrt{2 \operatorname{tr}\left( \mathbf{W} \boldsymbol{\Sigma}_1 \right)^2  }  \sqrt{\sum_{i=1}^k \frac{(k-1)^2}{n_i\left(n_i-1\right)}  +
	\sum_{i\neq l}^k \frac{1}{n_i n_l}  } } \\
\geq& \frac{\alpha^{2} p^{2 \delta} \nu^{2}+\nu^{2} \sum_{i=1}^{p^{\delta} } \omega_{i}^{2}}{\lambda_{p}^{*}\sqrt{2\left(\sum_{i=2}^{p-1} \omega_{i}^{4}+2\omega_{p}^4+2p \omega_{p}^{2} \alpha^{2} + p^{2} \alpha^{4} \right)}   \sqrt{\sum_{i=1}^k \frac{(k-1)^2}{n_i\left(n_i-1\right)}  +
	\sum_{i\neq l}^k \frac{1}{n_i n_l}  }  } .
\end{align*}
Based on this lower bound, we can obtain the following corollary:
\begin{cor}\label{special case2: Power of Corollary of CLT of main results}
Assume $\boldsymbol{\Sigma}_{1} = \boldsymbol{\Sigma}_{2}= \cdots=\boldsymbol{\Sigma}_{k}$, $\alpha_{1}=\cdots=\alpha_{p}=\alpha$,  $0 < \inf_{i}\omega_{i} \leq \sup_{i}\omega_i < \infty$, $\alpha^{2}=O\left(p^{-3 / 4}\right)$, $\lambda_{p}^{*} \sqrt{\sum_{i=1}^k \frac{(k-1)^2}{n_i\left(n_i-1\right)}  +
\sum_{i\neq l}^k \frac{1}{n_i n_l}  } =o\left( \nu^{2} p^{2 \delta -1}\right)$ and Assumptions A-D hold. Then under the alternative hypothesis \eqref{alternative}, we have
\begin{align*}
\lim _{ n, p \rightarrow \infty} \mathrm{P}\left( T_{n } \geq  \hat{\sigma}_{n, k} z_{\vartheta} \right) 
= 1.
\end{align*}
\end{cor}
From Corollary \ref{special case2: Power of Corollary of CLT of main results}, it is a startling finding that the power of the test asymptotically reaches 1 when $\delta > 1/2$, if all of $\boldsymbol{\Sigma}_{1}$'s eigenvalues are restricted away from $0$ and $\nu^{2} (\sum_{i=1}^k \frac{(k-1)^2}{n_i\left(n_i-1\right)}  +
\sum_{i\neq l}^k \frac{1}{n_i n_l})^{-1/2}  = O \left( 1 \right)$. In other words, extremely sparse signals with $\delta <1/2$ may not be handled well by the newly proposed test unless the sparse nonzero signals are very strong.

\section{Simulation Studies}\label{section3}
In this section, we conduct simulation studies to illustrate the performance of the proposed test (referred to as the $T_{w}$ test in the following context), compared to the tests proposed by \cite{hu2017on} (referred to as the $T_{hb}$ test), \cite{he2023high} (referred to as the $T_{hs}$ test). Similar to \cite{jiang022nonparametric}, we set 
$\alpha_{1}=\cdots=\alpha_{p}=\sqrt{5}p^{-3/8}$  and  $\omega_{q}  = \sqrt{2} (1 + \frac{q}{2p})$, $q=1, \ldots , p$, for convenience. In the experiments, the nominal significance level is set to $\vartheta=0.05$. For each fixed parameter setup, $5000$ replications are simulated to calculate empirical size and power.


The number of population is set to be $k =3$, and several different dimensions and sample sizes are considered in simulation. We set $n_{1}= 0.8n^{*}$, $n_{2}= n^{*}$ and $n_{3}= 1.2n^{*}$ with $n^{*} = 60, 100$, respectively, and the dimension $p=200,500,800$. Let the random samples are generated from the following model:
\begin{align}
\mathbf{x}_{i j}=\boldsymbol{\mu}_i +\boldsymbol{\Sigma}_{i}^{1/2} \mathbf{z}_{i j}, i \in\{1, 2, 3\}, j= 1, \ldots, n_i,
\end{align}
where $\left\{\mathbf{z}_{ij}\right\}_{j=1}^{n_i}$ are independent and identically distributed (i.i.d.) vector with coordinate $z_{i  j \theta }, \theta  \in \{1,\ldots, m \}$. We consider the following three distributions for $z_{i  j \theta }$.
\begin{itemize}
\item  The standard normal  $\mathcal{N}(0,1)$;

\item   The standardized Chi-squared distribution with degrees of freedom 2 , that is,  $2^{-1}\left\{\chi^{2}(2)-2\right\}$;

\item   The standardized  t-distribution with degrees of freedom 4 , that is,  $(2)^{-1 / 2} t(4)$.

\end{itemize}
We assigned  $\boldsymbol{\mu}_{1}=\boldsymbol{\mu}_{2}=\boldsymbol{\mu}_{3}=\mathbf{0}$  under $\mathrm{H}_{0}$ and also consider a possible alternative under $\mathrm{H}_{1}$,  $\boldsymbol{\mu}_{2}=\boldsymbol{\mu}_{3}=\mathbf{0}$, and $\boldsymbol{\mu}_{1}= (\overbrace{\tau , \ldots, \tau }^{ [p^{1-\rho}] }, 0, \ldots, 0 )^{\top}$, where $[a]$ denotes the integer part of $a$. In this setting, $\rho \in \left[0,1\right]$ controls the sparsity of signal, while the magnitude of $\tau $ gives the strength of signal in the experiments. We select the nonzero entries $\tau $ to be $\sqrt{2 r(1 / n_{1}+1 / n_{2} + 1 / n_{3}) \log p}$, where $r>0$. 

To gain perspectives on the level of sparsity in the simulation, we note that for $p = 500$ with $\rho = 0.04$, there are $500^{0.06}\approx 41$ signals. Similarly, for $p = 800$, the calculation yields $ 800 ^{0.06} \approx  55$ signals, reflecting a moderately level of sparsity. Therefore, we consider the paremeter  $\rho=\{0.1, 0.2, 0.3, 0.4\}$, covering highly dense signals for an alternative hypothesis at $\rho=0.1$, to moderately dense signals at $\rho=0.2$ or $\rho=0.3$, and finally to moderately sparse signals at $\rho=0.4$. Meanwhile, we set the signal strength $r=\{0.02, 0.04, 0.06, 0.08\}$. The above settings in our simulation studies are similar to those used in \cite{xu2016adaptive}.

For the covariance matrix, we consider the following two scenarios:

\textit{Scenario 1}: Unequal covariance matrices,  $\boldsymbol{\Sigma}_{1}=\mathbf{I}_{p}$, $\boldsymbol{\Sigma}_{2}=\left(0.5^{|i-j|}\right)$  and  $\boldsymbol{\Sigma}_{3}=\boldsymbol{\Sigma}_{1}+\boldsymbol{\Sigma}_{2}$.

\textit{Scenario 2}: Common covariance matrices,  $\boldsymbol{\Sigma}_{1}=\boldsymbol{\Sigma}_{2}=\boldsymbol{\Sigma}_{3}= \boldsymbol{\Sigma}$, where  $\boldsymbol{\Sigma}= \left(0.5^{|i-j|}\right)$.


The empirical sizes and powers are shown in Table \ref{Empirical sizes for Scenario 1}- \ref{Empirical sizes for Scenario 2} and \ref{Empirical powers for Scenario 1 and case 1}	-\ref{Empirical powers for Scenario 2 and case 3}, respectively. These results reveal that:
\begin{enumerate}
\item $T_{w}, T_{hb}$ and $T_{hs}$ could control the size well around $0.05$ for all cases. Although the size of $T_w$ is slightly larger than 0.05, this deviation decreases with increasing $n$ and $p$.
\item Our proposed test has the greatest power with $\rho =0.01$ and $r=0.02$.  This indicates that the weighted $L_{2}$-norm-based test is effective when the nonzero signals of the difference between mean vectors are weakly dense.
\item With an increase in dimension $p$, the empirical powers of all methods grow. The simulation results are consistent with  Corollary \ref{special case2: Power of Corollary of CLT of main results}.
\item When the nonzero signals of mean vectors become sparse, all the test's empirical power drastically decreases. Yet, our proposed test still outperforms other existing tests in the majority of scenarios.
\end{enumerate}
\begin{table}[htpb]
\centering
\caption{Empirical sizes for \it{Scenario 1}.} 
\label{Empirical sizes for Scenario 1}
\vspace{0.0 em} \centering
\scalebox{0.69}{\begin{tabular}{ccccccccccc}
	\hline
	&&  \multicolumn{3}{c}{$z_{i j \theta } \stackrel{i . i . d}{\sim} \mathcal{N}(0,1)$} &  \multicolumn{3}{c}{$z_{i j \theta } \stackrel{i . i . d}{\sim}\left(\chi_1^2-2\right) / 2$} & \multicolumn{3}{c}{$z_{i j \theta } \stackrel{i . i . d}{\sim} \left(t_4 / \sqrt{2}\right)$}  \\
	\cmidrule(r){3-5} \cmidrule(r){6-8} \cmidrule(r){9-11} 
	$p$ &$n^{*}$  &$T_{w}$&$T_{hb}$&$T_{hs}$&$T_{w}$&$T_{hb}$&$T_{hs}$&$T_{w}$&$T_{hb}$&$T_{hs}$ \\
	\hline
	200 & 60  & 0.060 & 0.056 & 0.051 & 0.063 & 0.054 & 0.045 & 0.060 & 0.058 & 0.043 \\
	& 100 & 0.059 & 0.062 & 0.055 & 0.060 & 0.059 & 0.051 & 0.058 & 0.052 & 0.044 \\
	500 & 60  & 0.057 & 0.054 & 0.052 & 0.057 & 0.050 & 0.045 & 0.059 & 0.058 & 0.045 \\
	& 100 & 0.052 & 0.050 & 0.045 & 0.061 & 0.053 & 0.050 & 0.059 & 0.053 & 0.048 \\
	800 & 60  & 0.056 & 0.051 & 0.052 & 0.056 & 0.058 & 0.053 & 0.059 & 0.054 & 0.044 \\
	& 100 & 0.055 & 0.048 & 0.052 & 0.059 & 0.051 & 0.046 & 0.055 & 0.050 & 0.048 \\
	\hline
\end{tabular}}
\end{table}

\begin{table}[htpb]
\centering
\caption{Empirical sizes for \it{Scenario 2}.} 
\label{Empirical sizes for Scenario 2}
\vspace{0.0 em} \centering
\scalebox{0.69}{\begin{tabular}{ccccccccccc}
	\hline
	&&  \multicolumn{3}{c}{$z_{i j \theta } \stackrel{i . i . d}{\sim} \mathcal{N}(0,1)$} &  \multicolumn{3}{c}{$z_{i j \theta } \stackrel{i . i . d}{\sim}\left(\chi_1^2-2\right) / 2$} & \multicolumn{3}{c}{$z_{i j \theta } \stackrel{i . i . d}{\sim} \left(t_4 / \sqrt{2}\right)$}  \\
	\cmidrule(r){3-5} \cmidrule(r){6-8} \cmidrule(r){9-11} 
	$p$ &$n^{*}$  &$T_{w}$&$T_{hb}$&$T_{hs}$&$T_{w}$&$T_{hb}$&$T_{hs}$&$T_{w}$&$T_{hb}$&$T_{hs}$ \\
	\hline
	200 & 60  & 0.064 & 0.057 & 0.056 & 0.069 & 0.060 & 0.052 & 0.066 & 0.057 & 0.050 \\
	& 100 & 0.059 & 0.055 & 0.052 & 0.061 & 0.058 & 0.051 & 0.062 & 0.059 & 0.050 \\
	500 & 60  & 0.058 & 0.049 & 0.050 & 0.065 & 0.056 & 0.052 & 0.066 & 0.059 & 0.051 \\
	& 100 & 0.058 & 0.050 & 0.048 & 0.058 & 0.055 & 0.050 & 0.060 & 0.052 & 0.049 \\
	800 & 60  & 0.054 & 0.049 & 0.049 & 0.061 & 0.056 & 0.050 & 0.058 & 0.052 & 0.047 \\
	& 100 & 0.056 & 0.049 & 0.054 & 0.057 & 0.056 & 0.048 & 0.059 & 0.054 & 0.045 \\
	\hline
\end{tabular}}
\end{table}

\begin{table}[htpb]
\centering
\caption{Empirical powers for \it{Scenario 1} when $z_{i j \theta } \stackrel{i . i . d}{\sim} \mathcal{N}(0,1)$.} 
\label{Empirical powers for Scenario 1 and case 1}
\vspace{0.0 em} \centering
\scalebox{0.69}{\begin{tabular}{ccccccccccccccc}
	\hline
	&& & \multicolumn{3}{c}{$r =0.02$} &  \multicolumn{3}{c}{$r=0.04$} & \multicolumn{3}{c}{$r=0.06$} & \multicolumn{3}{c}{$r=0.08$}  \\
	\cmidrule(r){4-6} \cmidrule(r){7-9} \cmidrule(r){10-12} \cmidrule(r){13-15}
	$p$ & $n^{*}$  & $\rho$  &$T_{w}$&$T_{hb}$&$T_{hs}$&$T_{w}$&$T_{hb}$&$T_{hs}$&$T_{w}$&$T_{hb}$&$T_{hs}$ &$T_{w}$&$T_{hb}$&$T_{hs}$ \\
	\hline
	200 & 60  & 0.1 & 0.874 & 0.357 & 0.305 & 0.998 & 0.746 & 0.638 & 1.000 & 0.944 & 0.873 & 1.000 & 0.994 & 0.970 \\
	&     & 0.2 & 0.420 & 0.191 & 0.171 & 0.765 & 0.431 & 0.373 & 0.930 & 0.680 & 0.579 & 0.988 & 0.856 & 0.748 \\
	&     & 0.3 & 0.176 & 0.124 & 0.110 & 0.328 & 0.230 & 0.201 & 0.496 & 0.366 & 0.316 & 0.656 & 0.507 & 0.430 \\
	&     & 0.4 & 0.098 & 0.085 & 0.081 & 0.154 & 0.141 & 0.131 & 0.209 & 0.197 & 0.176 & 0.281 & 0.272 & 0.237 \\
	500 & 60  & 0.1 & 0.999 & 0.671 & 0.559 & 1.000 & 0.987 & 0.952 & 1.000 & 1.000 & 0.999 & 1.000 & 1.000 & 1.000 \\
	&     & 0.2 & 0.720 & 0.312 & 0.268 & 0.977 & 0.711 & 0.615 & 0.999 & 0.930 & 0.851 & 1.000 & 0.992 & 0.970 \\
	&     & 0.3 & 0.259 & 0.167 & 0.148 & 0.512 & 0.332 & 0.294 & 0.743 & 0.565 & 0.490 & 0.888 & 0.758 & 0.660 \\
	&     & 0.4 & 0.112 & 0.103 & 0.093 & 0.178 & 0.166 & 0.151 & 0.268 & 0.251 & 0.224 & 0.376 & 0.383 & 0.319 \\
	800 & 60  & 0.1 & 1.000 & 0.836 & 0.748 & 1.000 & 0.999 & 0.994 & 1.000 & 1.000 & 1.000 & 1.000 & 1.000 & 1.000 \\
	&     & 0.2 & 0.859 & 0.404 & 0.352 & 0.998 & 0.866 & 0.769 & 1.000 & 0.990 & 0.960 & 1.000 & 1.000 & 0.996 \\
	&     & 0.3 & 0.287 & 0.185 & 0.164 & 0.623 & 0.416 & 0.343 & 0.843 & 0.669 & 0.585 & 0.957 & 0.871 & 0.771 \\
	&     & 0.4 & 0.124 & 0.113 & 0.098 & 0.207 & 0.186 & 0.169 & 0.314 & 0.299 & 0.257 & 0.430 & 0.431 & 0.375 \\
	200 & 100 & 0.1 & 0.876 & 0.347 & 0.292 & 0.997 & 0.755 & 0.654 & 1.000 & 0.955 & 0.876 & 1.000 & 0.996 & 0.972 \\
	&     & 0.2 & 0.434 & 0.188 & 0.166 & 0.754 & 0.418 & 0.357 & 0.941 & 0.687 & 0.570 & 0.982 & 0.851 & 0.747 \\
	&     & 0.3 & 0.187 & 0.135 & 0.116 & 0.323 & 0.225 & 0.206 & 0.501 & 0.371 & 0.314 & 0.642 & 0.512 & 0.425 \\
	&     & 0.4 & 0.101 & 0.095 & 0.094 & 0.156 & 0.131 & 0.124 & 0.215 & 0.196 & 0.166 & 0.280 & 0.274 & 0.235 \\
	500 & 100 & 0.1 & 0.998 & 0.657 & 0.558 & 1.000 & 0.986 & 0.950 & 1.000 & 1.000 & 0.999 & 1.000 & 1.000 & 1.000 \\
	&     & 0.2 & 0.715 & 0.318 & 0.262 & 0.980 & 0.713 & 0.613 & 1.000 & 0.938 & 0.862 & 1.000 & 0.994 & 0.966 \\
	&     & 0.3 & 0.254 & 0.162 & 0.154 & 0.508 & 0.338 & 0.285 & 0.743 & 0.560 & 0.468 & 0.886 & 0.752 & 0.654 \\
	&     & 0.4 & 0.115 & 0.102 & 0.085 & 0.177 & 0.168 & 0.148 & 0.268 & 0.255 & 0.218 & 0.379 & 0.368 & 0.307 \\
	800 & 100 & 0.1 & 1.000 & 0.841 & 0.748 & 1.000 & 1.000 & 0.995 & 1.000 & 1.000 & 1.000 & 1.000 & 1.000 & 1.000 \\
	&     & 0.2 & 0.859 & 0.407 & 0.352 & 0.998 & 0.857 & 0.767 & 1.000 & 0.988 & 0.959 & 1.000 & 1.000 & 0.994 \\
	&     & 0.3 & 0.293 & 0.184 & 0.158 & 0.625 & 0.423 & 0.367 & 0.841 & 0.670 & 0.582 & 0.956 & 0.875 & 0.773 \\
	&     & 0.4 & 0.116 & 0.105 & 0.097 & 0.203 & 0.188 & 0.154 & 0.328 & 0.313 & 0.270 & 0.437 & 0.430 & 0.367 \\
	\hline
\end{tabular}
}
\end{table}

\begin{table*}[htpb]
\centering
\caption{Empirical powers for \it{Scenario 1} when $z_{i j \theta } \stackrel{i . i . d}{\sim}\left(\chi_1^2-2\right) / 2 $.} 
\vspace{1.0 em} \centering
\scalebox{0.69}{	\begin{tabular}{ccccccccccccccc}
	\hline
	&& & \multicolumn{3}{c}{$r =0.02$} &  \multicolumn{3}{c}{$r=0.04$} & \multicolumn{3}{c}{$r=0.06$} & \multicolumn{3}{c}{$r=0.08$}  \\
	\cmidrule(r){4-6} \cmidrule(r){7-9} \cmidrule(r){10-12} \cmidrule(r){13-15}
	$p$ & $n^{*}$  & $\rho$  &$T_{w}$&$T_{hb}$&$T_{hs}$&$T_{w}$&$T_{hb}$&$T_{hs}$&$T_{w}$&$T_{hb}$&$T_{hs}$ &$T_{w}$&$T_{hb}$&$T_{hs}$ \\
	\hline
	200 & 60  & 0.1 & 0.878 & 0.348 & 0.304 & 0.997 & 0.761 & 0.631 & 1.000 & 0.951 & 0.875 & 1.000 & 0.995 & 0.964 \\
	&     & 0.2 & 0.411 & 0.186 & 0.155 & 0.775 & 0.431 & 0.357 & 0.942 & 0.679 & 0.569 & 0.986 & 0.852 & 0.734 \\
	&     & 0.3 & 0.165 & 0.116 & 0.102 & 0.324 & 0.227 & 0.182 & 0.489 & 0.357 & 0.302 & 0.626 & 0.493 & 0.410 \\
	&     & 0.4 & 0.095 & 0.089 & 0.076 & 0.158 & 0.145 & 0.123 & 0.220 & 0.202 & 0.169 & 0.278 & 0.275 & 0.223 \\
	500 & 60  & 0.1 & 0.998 & 0.659 & 0.560 & 1.000 & 0.988 & 0.950 & 1.000 & 1.000 & 0.997 & 1.000 & 1.000 & 1.000 \\
	&     & 0.2 & 0.711 & 0.309 & 0.260 & 0.976 & 0.703 & 0.583 & 0.999 & 0.932 & 0.846 & 1.000 & 0.995 & 0.967 \\
	&     & 0.3 & 0.245 & 0.157 & 0.133 & 0.525 & 0.343 & 0.274 & 0.734 & 0.538 & 0.458 & 0.880 & 0.758 & 0.641 \\
	&     & 0.4 & 0.114 & 0.097 & 0.080 & 0.181 & 0.169 & 0.140 & 0.273 & 0.261 & 0.216 & 0.370 & 0.360 & 0.298 \\
	800 & 60  & 0.1 & 1.000 & 0.837 & 0.741 & 1.000 & 1.000 & 0.994 & 1.000 & 1.000 & 1.000 & 1.000 & 1.000 & 1.000 \\
	&     & 0.2 & 0.858 & 0.401 & 0.334 & 0.997 & 0.852 & 0.753 & 1.000 & 0.987 & 0.952 & 1.000 & 1.000 & 0.995 \\
	&     & 0.3 & 0.284 & 0.175 & 0.148 & 0.617 & 0.420 & 0.339 & 0.845 & 0.673 & 0.569 & 0.951 & 0.861 & 0.760 \\
	&     & 0.4 & 0.118 & 0.109 & 0.097 & 0.209 & 0.188 & 0.157 & 0.311 & 0.299 & 0.246 & 0.438 & 0.427 & 0.356 \\
	200 & 100 & 0.1 & 0.879 & 0.342 & 0.285 & 0.997 & 0.754 & 0.632 & 1.000 & 0.951 & 0.872 & 1.000 & 0.996 & 0.973 \\
	&     & 0.2 & 0.423 & 0.195 & 0.170 & 0.763 & 0.421 & 0.354 & 0.937 & 0.670 & 0.551 & 0.987 & 0.857 & 0.743 \\
	&     & 0.3 & 0.174 & 0.125 & 0.110 & 0.309 & 0.218 & 0.184 & 0.496 & 0.365 & 0.307 & 0.655 & 0.521 & 0.428 \\
	&     & 0.4 & 0.094 & 0.092 & 0.081 & 0.152 & 0.137 & 0.127 & 0.218 & 0.205 & 0.177 & 0.281 & 0.276 & 0.229 \\
	500 & 100 & 0.1 & 0.998 & 0.661 & 0.560 & 1.000 & 0.989 & 0.945 & 1.000 & 1.000 & 0.999 & 1.000 & 1.000 & 1.000 \\
	&     & 0.2 & 0.715 & 0.317 & 0.260 & 0.975 & 0.717 & 0.609 & 0.999 & 0.936 & 0.860 & 1.000 & 0.991 & 0.961 \\
	&     & 0.3 & 0.249 & 0.155 & 0.139 & 0.507 & 0.330 & 0.277 & 0.741 & 0.571 & 0.474 & 0.883 & 0.756 & 0.639 \\
	&     & 0.4 & 0.115 & 0.101 & 0.093 & 0.173 & 0.168 & 0.144 & 0.277 & 0.265 & 0.222 & 0.381 & 0.366 & 0.311 \\
	800 & 100 & 0.1 & 1.000 & 0.846 & 0.748 & 1.000 & 0.999 & 0.993 & 1.000 & 1.000 & 1.000 & 1.000 & 1.000 & 1.000 \\
	&     & 0.2 & 0.859 & 0.406 & 0.349 & 0.997 & 0.864 & 0.758 & 1.000 & 0.988 & 0.955 & 1.000 & 1.000 & 0.997 \\
	&     & 0.3 & 0.310 & 0.182 & 0.149 & 0.630 & 0.424 & 0.354 & 0.835 & 0.669 & 0.580 & 0.951 & 0.864 & 0.775 \\
	&     & 0.4 & 0.116 & 0.113 & 0.096 & 0.213 & 0.198 & 0.170 & 0.323 & 0.302 & 0.258 & 0.445 & 0.436 & 0.363 \\
	\hline
\end{tabular}
}
\end{table*}	

\begin{table*}[htpb]
\centering
\caption{Empirical powers for \it{Scenario 1} when $z_{i j \theta } \stackrel{i . i . d}{\sim}\left(t_4 / \sqrt{2}\right)$. } 
\vspace{0.0 em} \centering
\scalebox{0.69}{	\begin{tabular}{ccccccccccccccc}
	\hline
	&& & \multicolumn{3}{c}{$r =0.02$} &  \multicolumn{3}{c}{$r=0.04$} & \multicolumn{3}{c}{$r=0.06$} & \multicolumn{3}{c}{$r=0.08$}  \\
	\cmidrule(r){4-6} \cmidrule(r){7-9} \cmidrule(r){10-12} \cmidrule(r){13-15}
	$p$ & $n^{*}$  & $\rho$  &$T_{w}$&$T_{hb}$&$T_{hs}$&$T_{w}$&$T_{hb}$&$T_{hs}$&$T_{w}$&$T_{hb}$&$T_{hs}$ &$T_{w}$&$T_{hb}$&$T_{hs}$ \\
	\hline
	200 & 60  & 0.1 & 0.877 & 0.355 & 0.292 & 0.995 & 0.756 & 0.613 & 1.000 & 0.947 & 0.860 & 1.000 & 0.997 & 0.963 \\
	&     & 0.2 & 0.422 & 0.199 & 0.151 & 0.770 & 0.427 & 0.335 & 0.930 & 0.662 & 0.540 & 0.985 & 0.851 & 0.729 \\
	&     & 0.3 & 0.176 & 0.128 & 0.097 & 0.327 & 0.236 & 0.181 & 0.503 & 0.370 & 0.278 & 0.650 & 0.509 & 0.405 \\
	&     & 0.4 & 0.101 & 0.092 & 0.070 & 0.162 & 0.140 & 0.106 & 0.205 & 0.204 & 0.155 & 0.279 & 0.269 & 0.219 \\
	500 & 60  & 0.1 & 0.999 & 0.664 & 0.535 & 1.000 & 0.987 & 0.944 & 1.000 & 1.000 & 0.996 & 1.000 & 1.000 & 0.999 \\
	&     & 0.2 & 0.715 & 0.304 & 0.240 & 0.970 & 0.700 & 0.582 & 0.999 & 0.932 & 0.840 & 1.000 & 0.993 & 0.959 \\
	&     & 0.3 & 0.246 & 0.155 & 0.124 & 0.505 & 0.335 & 0.258 & 0.744 & 0.555 & 0.441 & 0.888 & 0.748 & 0.609 \\
	&     & 0.4 & 0.113 & 0.105 & 0.085 & 0.189 & 0.167 & 0.129 & 0.276 & 0.265 & 0.203 & 0.381 & 0.372 & 0.290 \\
	800 & 60  & 0.1 & 1.000 & 0.840 & 0.720 & 1.000 & 0.999 & 0.991 & 1.000 & 1.000 & 0.997 & 1.000 & 1.000 & 1.000 \\
	&     & 0.2 & 0.860 & 0.394 & 0.317 & 0.998 & 0.856 & 0.751 & 1.000 & 0.985 & 0.944 & 1.000 & 1.000 & 0.993 \\
	&     & 0.3 & 0.294 & 0.184 & 0.142 & 0.612 & 0.412 & 0.319 & 0.852 & 0.688 & 0.554 & 0.953 & 0.858 & 0.743 \\
	&     & 0.4 & 0.124 & 0.117 & 0.081 & 0.211 & 0.193 & 0.144 & 0.324 & 0.307 & 0.234 & 0.431 & 0.433 & 0.334 \\
	200 & 100 & 0.1 & 0.876 & 0.351 & 0.285 & 0.997 & 0.755 & 0.632 & 1.000 & 0.954 & 0.870 & 1.000 & 0.995 & 0.969 \\
	&     & 0.2 & 0.424 & 0.198 & 0.157 & 0.767 & 0.430 & 0.341 & 0.935 & 0.681 & 0.551 & 0.986 & 0.858 & 0.724 \\
	&     & 0.3 & 0.177 & 0.128 & 0.104 & 0.328 & 0.238 & 0.190 & 0.490 & 0.359 & 0.293 & 0.646 & 0.512 & 0.414 \\
	&     & 0.4 & 0.104 & 0.086 & 0.073 & 0.161 & 0.149 & 0.122 & 0.212 & 0.194 & 0.164 & 0.280 & 0.287 & 0.230 \\
	500 & 100 & 0.1 & 0.999 & 0.663 & 0.541 & 1.000 & 0.986 & 0.947 & 1.000 & 1.000 & 0.997 & 1.000 & 1.000 & 1.000 \\
	&     & 0.2 & 0.713 & 0.312 & 0.248 & 0.974 & 0.705 & 0.581 & 0.999 & 0.938 & 0.847 & 1.000 & 0.994 & 0.964 \\
	&     & 0.3 & 0.254 & 0.153 & 0.123 & 0.510 & 0.336 & 0.275 & 0.753 & 0.556 & 0.455 & 0.885 & 0.766 & 0.631 \\
	&     & 0.4 & 0.108 & 0.105 & 0.085 & 0.183 & 0.163 & 0.133 & 0.279 & 0.263 & 0.204 & 0.379 & 0.369 & 0.303 \\
	800 & 100 & 0.1 & 1.000 & 0.838 & 0.729 & 1.000 & 0.999 & 0.993 & 1.000 & 1.000 & 0.999 & 1.000 & 1.000 & 1.000 \\
	&     & 0.2 & 0.866 & 0.414 & 0.327 & 0.997 & 0.852 & 0.747 & 1.000 & 0.990 & 0.949 & 1.000 & 0.999 & 0.992 \\
	&     & 0.3 & 0.299 & 0.185 & 0.145 & 0.616 & 0.414 & 0.345 & 0.858 & 0.693 & 0.562 & 0.952 & 0.862 & 0.749 \\
	&     & 0.4 & 0.120 & 0.110 & 0.090 & 0.205 & 0.183 & 0.152 & 0.316 & 0.303 & 0.245 & 0.429 & 0.432 & 0.348 \\
	\hline
\end{tabular}
}
\end{table*}

\begin{table*}[htpb]
\centering
\caption{Empirical powers for \it{Scenario 2} when $z_{i j \theta } \stackrel{i . i . d}{\sim} \mathcal{N}(0,1)$.} 
\vspace{ 0.0 em} \centering
\scalebox{0.69}{	\begin{tabular}{ccccccccccccccc}
	\hline
	&& & \multicolumn{3}{c}{$r =0.02$} &  \multicolumn{3}{c}{$r=0.04$} & \multicolumn{3}{c}{$r=0.06$} & \multicolumn{3}{c}{$r=0.08$}  \\
	\cmidrule(r){4-6} \cmidrule(r){7-9} \cmidrule(r){10-12} \cmidrule(r){13-15}
	$p$ & $n^{*}$  & $\rho$  &$T_{w}$&$T_{hb}$&$T_{hs}$&$T_{w}$&$T_{hb}$&$T_{hs}$&$T_{w}$&$T_{hb}$&$T_{hs}$ &$T_{w}$&$T_{hb}$&$T_{hs}$ \\
	\hline
	200 & 60  & 0.1 & 0.803 & 0.381 & 0.355 & 0.979 & 0.755 & 0.702 & 0.999 & 0.946 & 0.917 & 1.000 & 0.993 & 0.980 \\
	&     & 0.2 & 0.386 & 0.227 & 0.206 & 0.678 & 0.454 & 0.420 & 0.867 & 0.703 & 0.643 & 0.951 & 0.846 & 0.796 \\
	&     & 0.3 & 0.173 & 0.128 & 0.125 & 0.316 & 0.255 & 0.244 & 0.451 & 0.406 & 0.373 & 0.585 & 0.540 & 0.493 \\
	&     & 0.4 & 0.108 & 0.099 & 0.091 & 0.157 & 0.158 & 0.147 & 0.215 & 0.235 & 0.207 & 0.264 & 0.315 & 0.279 \\
	500 & 60  & 0.1 & 0.993 & 0.690 & 0.636 & 1.000 & 0.983 & 0.970 & 1.000 & 1.000 & 1.000 & 1.000 & 1.000 & 1.000 \\
	&     & 0.2 & 0.642 & 0.351 & 0.325 & 0.935 & 0.738 & 0.696 & 0.994 & 0.941 & 0.914 & 0.999 & 0.988 & 0.982 \\
	&     & 0.3 & 0.248 & 0.171 & 0.167 & 0.475 & 0.382 & 0.352 & 0.675 & 0.593 & 0.547 & 0.823 & 0.771 & 0.734 \\
	&     & 0.4 & 0.117 & 0.113 & 0.104 & 0.181 & 0.188 & 0.174 & 0.263 & 0.299 & 0.263 & 0.373 & 0.410 & 0.383 \\
	800 & 60  & 0.1 & 1.000 & 0.861 & 0.815 & 1.000 & 0.999 & 0.997 & 1.000 & 1.000 & 1.000 & 1.000 & 1.000 & 1.000 \\
	&     & 0.2 & 0.790 & 0.441 & 0.420 & 0.989 & 0.872 & 0.829 & 1.000 & 0.989 & 0.977 & 1.000 & 0.999 & 0.998 \\
	&     & 0.3 & 0.312 & 0.215 & 0.195 & 0.576 & 0.453 & 0.429 & 0.791 & 0.712 & 0.667 & 0.913 & 0.890 & 0.844 \\
	&     & 0.4 & 0.127 & 0.114 & 0.111 & 0.204 & 0.205 & 0.195 & 0.304 & 0.329 & 0.308 & 0.418 & 0.478 & 0.444 \\
	200 & 100 & 0.1 & 0.795 & 0.386 & 0.354 & 0.980 & 0.765 & 0.713 & 1.000 & 0.944 & 0.905 & 1.000 & 0.991 & 0.981 \\
	&     & 0.2 & 0.380 & 0.218 & 0.192 & 0.690 & 0.464 & 0.419 & 0.869 & 0.702 & 0.641 & 0.952 & 0.855 & 0.798 \\
	&     & 0.3 & 0.179 & 0.139 & 0.129 & 0.311 & 0.261 & 0.233 & 0.448 & 0.411 & 0.366 & 0.580 & 0.545 & 0.492 \\
	&     & 0.4 & 0.107 & 0.099 & 0.099 & 0.154 & 0.170 & 0.148 & 0.210 & 0.242 & 0.218 & 0.258 & 0.306 & 0.284 \\
	500 & 100 & 0.1 & 0.993 & 0.697 & 0.652 & 1.000 & 0.985 & 0.974 & 1.000 & 1.000 & 0.999 & 1.000 & 1.000 & 1.000 \\
	&     & 0.2 & 0.649 & 0.343 & 0.328 & 0.933 & 0.722 & 0.680 & 0.994 & 0.936 & 0.910 & 0.999 & 0.989 & 0.983 \\
	&     & 0.3 & 0.256 & 0.183 & 0.167 & 0.462 & 0.372 & 0.342 & 0.670 & 0.592 & 0.552 & 0.815 & 0.775 & 0.724 \\
	&     & 0.4 & 0.114 & 0.110 & 0.101 & 0.184 & 0.189 & 0.178 & 0.269 & 0.298 & 0.275 & 0.351 & 0.402 & 0.382 \\
	800 & 100 & 0.1 & 1.000 & 0.862 & 0.826 & 1.000 & 1.000 & 0.998 & 1.000 & 1.000 & 1.000 & 1.000 & 1.000 & 1.000 \\
	&     & 0.2 & 0.797 & 0.441 & 0.411 & 0.986 & 0.871 & 0.837 & 1.000 & 0.987 & 0.978 & 1.000 & 0.999 & 0.998 \\
	&     & 0.3 & 0.292 & 0.199 & 0.182 & 0.578 & 0.456 & 0.431 & 0.784 & 0.706 & 0.672 & 0.918 & 0.884 & 0.841 \\
	&     & 0.4 & 0.122 & 0.113 & 0.106 & 0.206 & 0.206 & 0.186 & 0.322 & 0.349 & 0.318 & 0.416 & 0.472 & 0.440 \\
	\hline
\end{tabular}
}
\end{table*}

\begin{table*}[htpb]
\centering
\caption{Empirical powers for \it{Scenario 2} when $z_{i j \theta } \stackrel{i . i . d}{\sim} \left(\chi_1^2-2\right) / 2$.} 
\vspace{0.0 em} \centering
\scalebox{0.69}{	\begin{tabular}{ccccccccccccccc}
	\hline
	&& & \multicolumn{3}{c}{$r =0.02$} &  \multicolumn{3}{c}{$r=0.04$} & \multicolumn{3}{c}{$r=0.06$} & \multicolumn{3}{c}{$r=0.08$}  \\
	\cmidrule(r){4-6} \cmidrule(r){7-9} \cmidrule(r){10-12} \cmidrule(r){13-15}
	$p$ & $n^{*}$  & $\rho$  &$T_{w}$&$T_{hb}$&$T_{hs}$&$T_{w}$&$T_{hb}$&$T_{hs}$&$T_{w}$&$T_{hb}$&$T_{hs}$ &$T_{w}$&$T_{hb}$&$T_{hs}$ \\
	\hline
	200 & 60  & 0.1 & 0.782 & 0.383 & 0.333 & 0.985 & 0.764 & 0.710 & 0.999 & 0.949 & 0.908 & 1.000 & 0.992 & 0.981 \\
	&     & 0.2 & 0.396 & 0.229 & 0.203 & 0.681 & 0.457 & 0.402 & 0.871 & 0.694 & 0.620 & 0.953 & 0.861 & 0.798 \\
	&     & 0.3 & 0.174 & 0.134 & 0.124 & 0.303 & 0.253 & 0.226 & 0.447 & 0.390 & 0.343 & 0.580 & 0.534 & 0.487 \\
	&     & 0.4 & 0.100 & 0.096 & 0.088 & 0.149 & 0.152 & 0.138 & 0.207 & 0.233 & 0.197 & 0.264 & 0.311 & 0.273 \\
	500 & 60  & 0.1 & 0.992 & 0.680 & 0.630 & 1.000 & 0.986 & 0.974 & 1.000 & 1.000 & 0.998 & 1.000 & 1.000 & 1.000 \\
	&     & 0.2 & 0.648 & 0.347 & 0.315 & 0.940 & 0.729 & 0.678 & 0.991 & 0.938 & 0.904 & 0.999 & 0.994 & 0.982 \\
	&     & 0.3 & 0.246 & 0.164 & 0.144 & 0.468 & 0.366 & 0.336 & 0.666 & 0.578 & 0.533 & 0.829 & 0.783 & 0.731 \\
	&     & 0.4 & 0.114 & 0.104 & 0.100 & 0.184 & 0.183 & 0.163 & 0.267 & 0.291 & 0.268 & 0.350 & 0.396 & 0.361 \\
	800 & 60  & 0.1 & 1.000 & 0.854 & 0.818 & 1.000 & 1.000 & 0.998 & 1.000 & 1.000 & 1.000 & 1.000 & 1.000 & 1.000 \\
	&     & 0.2 & 0.790 & 0.445 & 0.404 & 0.986 & 0.868 & 0.824 & 1.000 & 0.988 & 0.974 & 1.000 & 1.000 & 0.999 \\
	&     & 0.3 & 0.277 & 0.192 & 0.178 & 0.576 & 0.455 & 0.412 & 0.776 & 0.707 & 0.663 & 0.910 & 0.879 & 0.839 \\
	&     & 0.4 & 0.122 & 0.122 & 0.112 & 0.204 & 0.198 & 0.196 & 0.316 & 0.341 & 0.309 & 0.412 & 0.463 & 0.425 \\
	200 & 100 & 0.1 & 0.797 & 0.388 & 0.338 & 0.984 & 0.763 & 0.708 & 0.999 & 0.941 & 0.904 & 1.000 & 0.993 & 0.981 \\
	&     & 0.2 & 0.367 & 0.204 & 0.195 & 0.688 & 0.467 & 0.412 & 0.871 & 0.692 & 0.635 & 0.953 & 0.854 & 0.795 \\
	&     & 0.3 & 0.180 & 0.139 & 0.128 & 0.301 & 0.256 & 0.235 & 0.445 & 0.397 & 0.360 & 0.567 & 0.540 & 0.485 \\
	&     & 0.4 & 0.105 & 0.102 & 0.091 & 0.152 & 0.162 & 0.145 & 0.211 & 0.229 & 0.210 & 0.270 & 0.309 & 0.273 \\
	500 & 100 & 0.1 & 0.991 & 0.695 & 0.645 & 1.000 & 0.986 & 0.969 & 1.000 & 1.000 & 0.999 & 1.000 & 1.000 & 1.000 \\
	&     & 0.2 & 0.645 & 0.346 & 0.313 & 0.942 & 0.748 & 0.694 & 0.994 & 0.945 & 0.910 & 1.000 & 0.991 & 0.980 \\
	&     & 0.3 & 0.244 & 0.167 & 0.157 & 0.470 & 0.365 & 0.324 & 0.679 & 0.610 & 0.565 & 0.820 & 0.782 & 0.726 \\
	&     & 0.4 & 0.120 & 0.111 & 0.105 & 0.173 & 0.191 & 0.170 & 0.271 & 0.295 & 0.270 & 0.355 & 0.406 & 0.372 \\
	800 & 100 & 0.1 & 1.000 & 0.868 & 0.820 & 1.000 & 0.999 & 0.998 & 1.000 & 1.000 & 1.000 & 1.000 & 1.000 & 1.000 \\
	&     & 0.2 & 0.796 & 0.450 & 0.410 & 0.987 & 0.876 & 0.837 & 1.000 & 0.986 & 0.975 & 1.000 & 0.999 & 0.999 \\
	&     & 0.3 & 0.300 & 0.202 & 0.180 & 0.587 & 0.463 & 0.422 & 0.777 & 0.705 & 0.663 & 0.907 & 0.879 & 0.843 \\
	&     & 0.4 & 0.120 & 0.118 & 0.104 & 0.218 & 0.221 & 0.196 & 0.319 & 0.335 & 0.314 & 0.420 & 0.476 & 0.439 \\
	\hline
\end{tabular}
}
\end{table*}	

\begin{table}[htpb]
\centering
\caption{Empirical powers for \it{Scenario 2} when $z_{i j \theta } \stackrel{i . i . d}{\sim}\left(t_4 / \sqrt{2}\right)$.} 
\label{Empirical powers for Scenario 2 and case 3}
\vspace{0.0 em} \centering
\scalebox{0.69}{	\begin{tabular}{ccccccccccccccc}
	\hline
	&& & \multicolumn{3}{c}{$r =0.02$} &  \multicolumn{3}{c}{$r=0.04$} & \multicolumn{3}{c}{$r=0.06$} & \multicolumn{3}{c}{$r=0.08$}  \\
	\cmidrule(r){4-6} \cmidrule(r){7-9} \cmidrule(r){10-12} \cmidrule(r){13-15}
	$p$ & $n^{*}$  & $\rho$  &$T_{w}$&$T_{hb}$&$T_{hs}$&$T_{w}$&$T_{hb}$&$T_{hs}$&$T_{w}$&$T_{hb}$&$T_{hs}$ &$T_{w}$&$T_{hb}$&$T_{hs}$ \\
	\hline
	200 & 60  & 0.1 & 0.789 & 0.389 & 0.334 & 0.983 & 0.770 & 0.701 & 0.999 & 0.947 & 0.901 & 1.000 & 0.991 & 0.978 \\
	&     & 0.2 & 0.391 & 0.226 & 0.197 & 0.682 & 0.453 & 0.401 & 0.870 & 0.692 & 0.611 & 0.948 & 0.845 & 0.776 \\
	&     & 0.3 & 0.172 & 0.134 & 0.114 & 0.315 & 0.257 & 0.225 & 0.455 & 0.402 & 0.345 & 0.578 & 0.547 & 0.481 \\
	&     & 0.4 & 0.096 & 0.095 & 0.083 & 0.156 & 0.151 & 0.131 & 0.201 & 0.230 & 0.205 & 0.274 & 0.319 & 0.272 \\
	500 & 60  & 0.1 & 0.992 & 0.689 & 0.620 & 1.000 & 0.982 & 0.961 & 1.000 & 1.000 & 0.998 & 1.000 & 1.000 & 0.999 \\
	&     & 0.2 & 0.642 & 0.341 & 0.299 & 0.940 & 0.731 & 0.674 & 0.993 & 0.938 & 0.903 & 0.999 & 0.989 & 0.976 \\
	&     & 0.3 & 0.247 & 0.175 & 0.151 & 0.466 & 0.378 & 0.334 & 0.676 & 0.593 & 0.520 & 0.823 & 0.781 & 0.705 \\
	&     & 0.4 & 0.119 & 0.114 & 0.096 & 0.191 & 0.184 & 0.166 & 0.251 & 0.289 & 0.255 & 0.353 & 0.404 & 0.344 \\
	800 & 60  & 0.1 & 1.000 & 0.859 & 0.803 & 1.000 & 0.999 & 0.996 & 1.000 & 1.000 & 1.000 & 1.000 & 1.000 & 1.000 \\
	&     & 0.2 & 0.793 & 0.438 & 0.387 & 0.984 & 0.863 & 0.810 & 1.000 & 0.990 & 0.970 & 1.000 & 1.000 & 0.997 \\
	&     & 0.3 & 0.293 & 0.209 & 0.175 & 0.579 & 0.460 & 0.409 & 0.783 & 0.714 & 0.639 & 0.913 & 0.880 & 0.830 \\
	&     & 0.4 & 0.122 & 0.110 & 0.096 & 0.205 & 0.215 & 0.181 & 0.300 & 0.327 & 0.287 & 0.417 & 0.476 & 0.411 \\
	200 & 100 & 0.1 & 0.799 & 0.390 & 0.339 & 0.982 & 0.769 & 0.704 & 0.998 & 0.949 & 0.900 & 1.000 & 0.990 & 0.976 \\
	&     & 0.2 & 0.380 & 0.222 & 0.186 & 0.686 & 0.457 & 0.410 & 0.867 & 0.698 & 0.612 & 0.951 & 0.862 & 0.794 \\
	&     & 0.3 & 0.164 & 0.129 & 0.116 & 0.307 & 0.253 & 0.227 & 0.443 & 0.402 & 0.347 & 0.571 & 0.553 & 0.489 \\
	&     & 0.4 & 0.103 & 0.100 & 0.088 & 0.151 & 0.167 & 0.139 & 0.214 & 0.234 & 0.211 & 0.267 & 0.305 & 0.267 \\
	500 & 100 & 0.1 & 0.993 & 0.690 & 0.627 & 1.000 & 0.986 & 0.970 & 1.000 & 0.999 & 0.999 & 1.000 & 1.000 & 1.000 \\
	&     & 0.2 & 0.642 & 0.347 & 0.303 & 0.945 & 0.738 & 0.676 & 0.993 & 0.933 & 0.898 & 0.999 & 0.992 & 0.978 \\
	&     & 0.3 & 0.249 & 0.167 & 0.159 & 0.481 & 0.371 & 0.332 & 0.674 & 0.605 & 0.545 & 0.826 & 0.772 & 0.718 \\
	&     & 0.4 & 0.119 & 0.111 & 0.095 & 0.192 & 0.198 & 0.166 & 0.260 & 0.282 & 0.249 & 0.350 & 0.402 & 0.350 \\
	800 & 100 & 0.1 & 1.000 & 0.859 & 0.804 & 1.000 & 0.999 & 0.996 & 1.000 & 1.000 & 1.000 & 1.000 & 1.000 & 1.000 \\
	&     & 0.2 & 0.806 & 0.452 & 0.410 & 0.989 & 0.875 & 0.823 & 1.000 & 0.987 & 0.971 & 1.000 & 1.000 & 0.998 \\
	&     & 0.3 & 0.305 & 0.201 & 0.178 & 0.582 & 0.466 & 0.412 & 0.792 & 0.715 & 0.653 & 0.918 & 0.877 & 0.830 \\
	&     & 0.4 & 0.125 & 0.115 & 0.100 & 0.219 & 0.216 & 0.189 & 0.310 & 0.346 & 0.297 & 0.416 & 0.472 & 0.423 \\
	\hline
\end{tabular}
}
\end{table}

\section{Conclusions and discussions}
In this paper, we propose a weighted $L_{2}$-norm-based test to examine the equality of the mean vectors of several populations. Our proposed test provides a unified framework with adjustable settings. Compared with other existing tests, our proposed test has some advantages when the non-zero signal is weakly dense. In our simulation studies, we set 
$\alpha_{1}=\cdots=\alpha_{p}=\sqrt{5}p^{-3/8}$  and  $\omega_{q}  = \sqrt{2} (1 + \frac{q}{2p})$, $q=1, \ldots , p$ as the weight function. However, there are other choices for the weight function. An interesting topic is still to consider $\mathbf {B} = \operatorname{diag} \left\{\omega_{1}^{2}, \omega_{2}^{2}, \ldots, \omega_{p}^{2}\right\}$, but to set $\boldsymbol{\alpha} =\left( \alpha_{1}, \ldots, \alpha_{p}\right)^{\top}$ is random. At this point, the test statistic can be regarded as a weighted sum of one test statistic related to the $L_{2}$-norm-based test and another test statistic related to random projection.  This prospect is indeed intriguing for
future research and the interesting work will be continued in our subsequent paper.

\section{Acknowledgements}
Bai's research was supported by National Natural Science Foundation of China (No. 12271536, No. 12171198), and Natural Science Foundation of Jilin Province (No. 20210101147JC); Niu's research was supported by the Natural Science Foundation of Shandong Province (No. ZR2021QA077). The first author thanks Professor Yunlu Jiang for helpful discussion.


\appendix
\section{Proofs of main results}\label{Appendix}
\setcounter{lemma}{0}
\renewcommand{\thelemma}{\Alph{section}\arabic{lemma}}
\subsection{Proof of Theorem \ref{CLT of main results} } 
Write
\begin{align}
T_{n}  = T_{n1} - T_{n2},
\end{align}
where
\begin{align*}
T_{n1} =& (k-1) \sum_{i=1}^{k} \frac{1}{n_{i}\left(n_{i}-1\right)} \sum_{j_{1} \neq j_{2}} \mathbf{x}_{i j_{1}}^{\top}  \mathbf{W} \mathbf{x}_{i j_{2}} \\
T_{n2}=& \sum_{i<l}^{k} \frac{2}{n_{i} n_{l}} \sum_{j_{1}}^{n_{i}} \sum_{j_{2}}^{n_{l}} \mathbf{x}_{i j_{1}}^{\top}  \mathbf{W} \mathbf{x}_{l j_{2}}.
\end{align*}
It is very easy to show that
\begin{align*}
\mathbb{E} \left( T_{n} \right) =& \mathbb{E} \left(  T_{n1} \right) -  \mathbb{E} \left( T_{n2} \right) \\
=&  (k-1) \sum_{i=1}^{k} \frac{1}{n_{i}\left(n_{i}-1\right)} \sum_{j_{1} \neq j_{2}} \mathbb{E} \left(\mathbf{x}_{i j_{1}}^{\top}  \mathbf{W} \mathbf{x}_{i j_{2}} \right) \\
&- \sum_{i<l}^{k} \frac{2}{n_{i} n_{l}} \sum_{j_{1}, j_{2}} \left( \mathbb{E} \mathbf{x}_{i j_{1}}\right)^{\top} \mathbf{W} \left( \mathbb{E} \mathbf{x}_{l j_{2}} \right) \\
=& (k-1) \sum_{i=1}^{k} \boldsymbol{\mu}_{i}^{\top} \mathbf{W} \boldsymbol{\mu}_{i}
-2  \sum_{i<l}^{k} \boldsymbol{\mu}_{i}^{\top} \mathbf{W} \boldsymbol{\mu}_{l} \\
= & \sum_{i<l}^{k}\left( \boldsymbol{\mu}_{i}- \boldsymbol{\mu}_{l} \right)^{\top} \mathbf{W}\left( \boldsymbol{\mu}_{i}- \boldsymbol{\mu}_{l} \right) ,
\end{align*} 
and
\begin{align*}
\operatorname{Var} \left( T_{n} \right) =& \operatorname{Var} \left(  T_{n1} \right) +  \operatorname{Var} \left( T_{n2} \right)  - \operatorname{Cov} \left( T_{n1}, T_{n2} \right) - \operatorname{Cov} \left( T_{n2}, T_{n1} \right).
\end{align*}
\begin{align}
&\operatorname{Var} \left(  T_{n1} \right) \notag \\
=&  \sum_{i_1 , i_2} \frac{(k-1)^2}{n_{i_1}\left(n_{i_1 }-1\right) n_{i_2 }\left(n_{i_2}-1\right)} \sum_{j_1\neq j_2} \sum_{j_3\neq j_4} \mathbb{E} \left( \mathbf{x}_{i_1 j_1}\right)^{\top} \mathbf{W} \mathbf{x}_{i_1 j_2}\left(\mathbf{x}_{i_2 j_3}\right)^{\top} \mathbf{W} \mathbf{x}_{i_2 j_4}  \notag \\
&-( (k-1) \sum_{i=1}^{k} \boldsymbol{\mu}_{i}^{\top} \mathbf{W} \boldsymbol{\mu}_{i}  )^{2} \notag \\
=&  \sum_{i_1= 1}\frac{(k-1)^2}{n_{i_1}^{2}\left(n_{i_1 }-1\right)^{2}}  \sum_{j_1\neq j_2} \sum_{j_3\neq j_4} \mathbb{E} \left( \mathbf{x}_{i_1 j_1}\right)^{\top} \mathbf{W} \mathbf{x}_{i_1 j_2}\left(\mathbf{x}_{i_1 j_3}\right)^{\top} \mathbf{W} \mathbf{x}_{i_1 j_4}     \notag \\
&+ \sum_{i_1 \neq i_2} (k-1)^2 \boldsymbol{\mu}_{i_1}^{\top} \mathbf{W} \boldsymbol{\mu}_{i_1}\boldsymbol{\mu}_{i_2}^{\top} \mathbf{W} \boldsymbol{\mu}_{i_2}             
-( (k-1) \sum_{i=1}^{k} \boldsymbol{\mu}_{i}^{\top} \mathbf{W} \boldsymbol{\mu}_{i}  )^{2} \notag \\
=& \sum_{i_1= 1}\frac{(k-1)^2}{n_{i_1}^{2}\left(n_{i_1 }-1\right)^{2}}  \mathbb{E} \left(\sum_{j_1\neq j_2 \neq j_3 }  4( \mathbf{x}_{i_1 j_1})^{\top} \mathbf{W} \mathbf{x}_{i_1 j_2}(\mathbf{x}_{i_1 j_3})^{\top} \mathbf{W} \mathbf{x}_{i_1 j_2}  + 2\sum_{j_1\neq j_2} ( ( \mathbf{x}_{i_1 j_1})^{\top} \mathbf{W} \mathbf{x}_{i_1 j_2} )^{2} \notag \right. \\
& \left.  +  \sum_{j_1\neq j_2 \neq  j_3\neq j_4}  \left( \mathbf{x}_{i_1 j_1}\right)^{\top} \mathbf{W} \mathbf{x}_{i_1 j_2}\left(\mathbf{x}_{i_1 j_3}\right)^{\top} \mathbf{W} \mathbf{x}_{i_1 j_4} \right)  + \sum_{i_1 \neq i_2} (k-1)^2 \boldsymbol{\mu}_{i_1}^{\top} \mathbf{W} \boldsymbol{\mu}_{i_1}\boldsymbol{\mu}_{i_2}^{\top} \mathbf{W} \boldsymbol{\mu}_{i_2}      \notag   \\    
&-( (k-1) \sum_{i=1}^{k} \boldsymbol{\mu}_{i}^{\top} \mathbf{W} \boldsymbol{\mu}_{i}  )^{2} \notag \\
=& \sum_{i_1= 1}\frac{(k-1)^2}{n_{i_1}^{2}\left(n_{i_1 }-1\right)^{2}}  \left(   4n_{i_1}\left(n_{i_1 }-1\right)\left(n_{i_1 }-2\right)\boldsymbol{\mu}_{i_1}^{\top} \mathbf{W} \boldsymbol{\Sigma}_{i_1} \mathbf{W} \boldsymbol{\mu}_{i_1}   \notag \right. \\
& \left. +4n_{i_1}\left(n_{i_1 }-1\right)\boldsymbol{\mu}_{i_1}^{\top} \mathbf{W} \boldsymbol{\Sigma}_{i_1} \mathbf{W} \boldsymbol{\mu}_{i_1}  + 2n_{i_1}\left(n_{i_1 }-1\right)\operatorname{tr}\left( \mathbf{W} \boldsymbol{\Sigma}_{i_{1}}\right)^2 
\right)              \notag \\
=&   (k-1)^2 \sum_{i_1= 1} ( \frac{4}{n_{i_1} }   \boldsymbol{\mu}_{i_1}^{\top} \mathbf{W} \boldsymbol{\Sigma}_{i_1} \mathbf{W} \boldsymbol{\mu}_{i_1}  + \frac{2}{n_{i_1}\left(n_{i_1 }-1\right)}\operatorname{tr}( \mathbf{W} \boldsymbol{\Sigma}_{i_{1}})^2 
).    \label{varTn1}      
\end{align}
Similarly, we can obtain
\begin{align}
&\operatorname{Var} \left( T_{n2} \right) \notag \\
= &  \sum_{i_{1}<l_{1}}^{k} \sum_{i_{2}<l_{2}}^{k} \frac{4}{n_{i_{1}} n_{l_{1}} n_{i_{2}} n_{l_{2}} } \sum_{j_{1}, j_{2}}\sum_{j_{3}, j_{4}} \mathbb{E} \mathbf{x}_{i_{1} j_{1}}^{\top}  \mathbf{W} \mathbf{x}_{l_{1} j_{2}} \mathbf{x}_{i_{2} j_{3}}^{\top}  \mathbf{W} \mathbf{x}_{l_{2} j_{4}}      \notag\\
&-  (-2  \sum_{i<l}^{k} \boldsymbol{\mu}_{i}^{\top} \mathbf{W} \boldsymbol{\mu}_{l})^{2}    \notag\\
=&    \sum_{i_1\neq l_1 \neq i_2 } \frac{4}{n_{i_{1}}^{2} n_{l_{1}} n_{i_{2}}  } \sum_{j_{1}, j_{2}}\sum_{j_{3}, j_{4}} \mathbb{E} \mathbf{x}_{i_{1} j_{1}}^{\top}  \mathbf{W} \mathbf{x}_{l_{1} j_{2}} \mathbf{x}_{i_{2} j_{3}}^{\top}  \mathbf{W} \mathbf{x}_{i_{1} j_{4}} - (-2  \sum_{i<l}^{k} \boldsymbol{\mu}_{i}^{\top} \mathbf{W} \boldsymbol{\mu}_{l} )^{2}  \notag\\
&+  \sum_{i_1\neq l_1 \neq  i_2\neq l_2} \frac{1}{n_{i_{1}} n_{l_{1}} n_{i_{2}} n_{l_{2}} } \sum_{j_{1}, j_{2}}\sum_{j_{3}, j_{4}} \mathbb{E} \left( \mathbf{x}_{i_1 j_1}\right)^{\top} \mathbf{W} \mathbf{x}_{l_1 j_2}\left(\mathbf{x}_{i_2 j_3}\right)^{\top} \mathbf{W} \mathbf{x}_{l_2 j_4}  \notag \\
& + \sum_{i_1\neq l_1} \frac{2}{n_{i_{1}}^{2} n_{l_{1}}^{2}  } \sum_{j_{1}, j_{2}}\sum_{j_{3}, j_{4}} \notag   \mathbb{E} \mathbf{x}_{i_{1} j_{1}}^{\top}  \mathbf{W} \mathbf{x}_{l_{1} j_{2}} \mathbf{x}_{i_{1} j_{3}}^{\top}  \mathbf{W} \mathbf{x}_{l_{1} j_{4}}  \notag \\
=&   \sum_{i_1\neq l_1 \neq i_2 } 4(  \frac{1}{n_{i_{1}}  } \boldsymbol{\mu}_{i_2}^{\top} \mathbf{W} \boldsymbol{\Sigma}_{i_1} \mathbf{W} \boldsymbol{\mu}_{l_1}+  \boldsymbol{\mu}_{i_1}^{\top} \mathbf{W}\boldsymbol{\mu}_{l_1} \boldsymbol{\mu}_{i_2}^{\top} \mathbf{W}\boldsymbol{\mu}_{i_1} )  + \sum_{i_1\neq l_1 \neq  i_2\neq l_2}\boldsymbol{\mu}_{i_1}^{\top} \mathbf{W}\boldsymbol{\mu}_{l_1} \boldsymbol{\mu}_{i_2}^{\top} \mathbf{W}\boldsymbol{\mu}_{l_2}  \notag \\
&  + 2\sum_{i_1\neq l_1} (   \frac{1}{n_{l_{1}}  } \boldsymbol{\mu}_{i_1}^{\top} \mathbf{W} \boldsymbol{\Sigma}_{l_1} \mathbf{W} \boldsymbol{\mu}_{i_1} + \frac{1}{n_{i_{1}}  }  \boldsymbol{\mu}_{l_1}^{\top} \mathbf{W} \boldsymbol{\Sigma}_{i_1} \mathbf{W} \boldsymbol{\mu}_{l_1} + \frac{1}{n_{i_{1}} n_{l_{1}}  } \operatorname{tr}( \mathbf{W} \boldsymbol{\Sigma}_{i_{1}} \mathbf{W} \boldsymbol{\Sigma}_{l_{1}}) )   \notag \\
&+ \sum_{i_1\neq l_1}2 \boldsymbol{\mu}_{i_{1}}^{\top} \mathbf{W} \boldsymbol{\mu}_{l_{1}}\boldsymbol{\mu}_{i_{1}}^{\top} \mathbf{W} \boldsymbol{\mu}_{l_{1}}  - (-2  \sum_{i<l}^{k} \boldsymbol{\mu}_{i}^{\top} \mathbf{W} \boldsymbol{\mu}_{l} )^{2}    \notag \\
=& \sum_{i_1\neq l_1} 2( \frac{1}{n_{l_{1}}  } \boldsymbol{\mu}_{i_1}^{\top} \mathbf{W} \boldsymbol{\Sigma}_{l_1} \mathbf{W} \boldsymbol{\mu}_{i_1} + \frac{1}{n_{i_{1}}  }  \boldsymbol{\mu}_{l_1}^{\top} \mathbf{W} \boldsymbol{\Sigma}_{i_1} \mathbf{W} \boldsymbol{\mu}_{l_1} + \frac{1}{n_{i_{1}} n_{l_{1}}  } \operatorname{tr}\left( \mathbf{W} \boldsymbol{\Sigma}_{i_{1}} \mathbf{W} \boldsymbol{\Sigma}_{l_{1}}\right) )     \notag  \\
& +  \sum_{i_1\neq l_1 \neq i_2 }   \frac{4}{n_{i_{1}}  } \boldsymbol{\mu}_{i_2}^{\top} \mathbf{W} \boldsymbol{\Sigma}_{i_1} \mathbf{W} \boldsymbol{\mu}_{l_1}     \notag \\
=&     \sum_{i_1\neq l_1} ( \frac{4}{n_{i_{1}}  }  \boldsymbol{\mu}_{l_1}^{\top} \mathbf{W} \boldsymbol{\Sigma}_{i_1} \mathbf{W} \boldsymbol{\mu}_{l_1}  + \frac{2}{n_{i_{1}} n_{l_{1}}  } \operatorname{tr}( \mathbf{W} \boldsymbol{\Sigma}_{i_{1}} \mathbf{W} \boldsymbol{\Sigma}_{l_{1}}) ) +      \sum_{i_1\neq l_1 \neq i_2 }   \frac{4}{n_{i_{1}}  } \boldsymbol{\mu}_{i_2}^{\top} \mathbf{W} \boldsymbol{\Sigma}_{i_1} \mathbf{W} \boldsymbol{\mu}_{l_1}.  \label{varTn2}
\end{align}
\begin{align}
\operatorname{Cov} \left( T_{n1}, T_{n2} \right)=	&
\sum_{i_{1}=1}^{k}\sum_{i_{2}<l_{2}}^{k} \frac{1}{n_{i_{1}}\left(n_{i_{1}}-1\right)} \frac{2 (k-1)}{ n_{i_{2}} n_{l_{2}} }  \sum_{j_{1} \neq j_{2}} \sum_{j_{3}, j_{4}} \mathbb{E} \mathbf{x}_{i_{1} j_{1}}^{\top}  \mathbf{W} \mathbf{x}_{i_{1} j_{2}}   \mathbf{x}_{i_{2} j_{3}}^{\top}  \mathbf{W} \mathbf{x}_{l_{2} j_{4}}   \notag \\
&- 2 (k-1) ( \sum_{i=1}^{k} \boldsymbol{\mu}_{i}^{\top} \mathbf{W} \boldsymbol{\mu}_{i}   ) ( \sum_{i<l}^{k} \boldsymbol{\mu}_{i}^{\top} \mathbf{W} \boldsymbol{\mu}_{l} )     \notag \\
=&   \sum_{i_{1}=i_{2} \neq l_{2}}^{k} \frac{1}{n_{i_{1}}^{2} \left(n_{i_{1}}-1\right)} \frac{ 2(k-1)}{  n_{l_{2}} } \sum_{j_{1} \neq j_{2}} \sum_{j_{3}, j_{4}} \mathbb{E} \left(  \mathbf{x}_{i_{1} j_{1}}^{\top}  \mathbf{W} \mathbf{x}_{i_{1} j_{2}}   \mathbf{x}_{i_{1} j_{3}}^{\top}  \mathbf{W} \mathbf{x}_{l_{2} j_{4}}   \right)       \notag \\
&+  \sum_{i_{1} \neq i_{2} \neq l_{2}}^{k} \frac{1}{n_{i_{1}}\left(n_{i_{1}}-1\right)} \frac{ (k-1)}{ n_{i_{2}} n_{l_{2}} }  \sum_{j_{1} \neq j_{2}} \sum_{j_{3}, j_{4}} \mathbb{E} \left( \mathbf{x}_{i_{1} j_{1}}^{\top}  \mathbf{W} \mathbf{x}_{i_{1} j_{2}}   \mathbf{x}_{i_{2} j_{3}}^{\top}  \mathbf{W} \mathbf{x}_{l_{2} j_{4}}    \right)         \notag \\
&- 2 (k-1) ( \sum_{i=1}^{k} \boldsymbol{\mu}_{i}^{\top} \mathbf{W} \boldsymbol{\mu}_{i}  ) ( \sum_{i<l}^{k} \boldsymbol{\mu}_{i}^{\top} \mathbf{W} \boldsymbol{\mu}_{l} )     \notag \\
=& \sum_{i_{1}=i_{2} \neq l_{2}}^{k} 2 (k-1)(  \boldsymbol{\mu}_{i_{1}}^{\top} \mathbf{W} \boldsymbol{\mu}_{i_{1}}\boldsymbol{\mu}_{i_{1}}^{\top} \mathbf{W} \boldsymbol{\mu}_{l_{2}} +  \frac{2}{n_{i_{1}} } \boldsymbol{\mu}_{i_1}^{\top} \mathbf{W} \boldsymbol{\Sigma}_{i_1} \mathbf{W} \boldsymbol{\mu}_{l_2}  ) \notag\\
& +  \sum_{i_{1} \neq i_{2} \neq l_{2}}^{k} \boldsymbol{\mu}_{i_{1}}^{\top} \mathbf{W} \boldsymbol{\mu}_{i_{1}}\boldsymbol{\mu}_{i_{2}}^{\top} \mathbf{W} \boldsymbol{\mu}_{l_{2}}    
- 2 (k-1) ( \sum_{i=1}^{k} \boldsymbol{\mu}_{i}^{\top} \mathbf{W} \boldsymbol{\mu}_{i}   ) ( \sum_{i<l}^{k} \boldsymbol{\mu}_{i}^{\top} \mathbf{W} \boldsymbol{\mu}_{l} )     \notag \\
=&   (k-1) \sum_{i_{1}=i_{2} \neq l_{2}}^{k} \frac{4}{n_{i_{1}} } \boldsymbol{\mu}_{i_1}^{\top} \mathbf{W} \boldsymbol{\Sigma}_{i_1} \mathbf{W} \boldsymbol{\mu}_{l_2}.         \label{covTn2Tn1}
\end{align}
Finally, from \eqref{varTn1}-\eqref{covTn2Tn1}, we have
\begin{align*}
\operatorname{Var}\left(T_{n}\right)= &    (k-1)^2 \sum_{i_1= 1}   \frac{4}{n_{i_1} }   \boldsymbol{\mu}_{i_1}^{\top} \mathbf{W} \boldsymbol{\Sigma}_{i_1} \mathbf{W} \boldsymbol{\mu}_{i_1}  + \frac{2}{n_{i_1}\left(n_{i_1 }-1\right)}\operatorname{tr}\left( \mathbf{W} \boldsymbol{\Sigma}_{i_{1}})\right)^2 
\\
&+   \sum_{i_1\neq l_1} \left( \frac{4}{n_{i_{1}}  }  \boldsymbol{\mu}_{l_1}^{\top} \mathbf{W} \boldsymbol{\Sigma}_{i_1} \mathbf{W} \boldsymbol{\mu}_{l_1}  + \frac{2}{n_{i_{1}} n_{l_{1}}  } \operatorname{tr}\left( \mathbf{W} \boldsymbol{\Sigma}_{i_{1}} \mathbf{W} \boldsymbol{\Sigma}_{l_{1}}\right) \right) \\
&+      \sum_{i_1\neq l_1 \neq i_2 }   \frac{4}{n_{i_{1}}  } \boldsymbol{\mu}_{i_2}^{\top} \mathbf{W} \boldsymbol{\Sigma}_{i_1} \mathbf{W} \boldsymbol{\mu}_{l_1}   -  8(k-1) \sum_{i_{1}=i_{2} \neq l_{2}}^{k} \frac{1}{n_{i_{1}} } \boldsymbol{\mu}_{i_1}^{\top} \mathbf{W} \boldsymbol{\Sigma}_{i_1} \mathbf{W} \boldsymbol{\mu}_{l_2}              \\	
= 
&   \sum_{i_1\neq l_1} \left( \frac{4}{n_{i_{1}}  }  \boldsymbol{\mu}_{l_1}^{\top} \mathbf{W} \boldsymbol{\Sigma}_{i_1} \mathbf{W} \boldsymbol{\mu}_{l_1}  \right) +   4   \sum_{i_1\neq l_1 \neq i_2 } \left(  \frac{1}{n_{i_{1}}  } \boldsymbol{\mu}_{i_2}^{\top} \mathbf{W} \boldsymbol{\Sigma}_{i_1} \mathbf{W} \boldsymbol{\mu}_{l_1}\right) \\
& +   \sum_{i_1= 1} \frac{4(k-1)^2}{n_{i_1} }   \boldsymbol{\mu}_{i_1}^{\top} \mathbf{W} \boldsymbol{\Sigma}_{i_1} \mathbf{W} \boldsymbol{\mu}_{i_1}- 
8(k-1) \sum_{i_{1}=i_{2} \neq l_{2}}^{k} \frac{1}{n_{i_{1}} } \boldsymbol{\mu}_{i_1}^{\top} \mathbf{W} \boldsymbol{\Sigma}_{i_1} \mathbf{W} \boldsymbol{\mu}_{l_2}  \\
&+  \sum_{i_1= 1}  \frac{2(k-1)^2}{n_{i_1}\left(n_{i_1 }-1\right)}\operatorname{tr}\left( \mathbf{W} \boldsymbol{\Sigma}_i\right)^2 
+   \sum_{i_1\neq l_1} \frac{2}{n_{i_{1}} n_{l_{1}}  } \operatorname{tr}\left( \mathbf{W} \boldsymbol{\Sigma}_{i_{1}} \mathbf{W} \boldsymbol{\Sigma}_{l_{1}}\right)  \\	
= &  \sum_{i=1}^k \frac{2(k-1)^2}{n_i\left(n_i-1\right)} \operatorname{tr}\left( \mathbf{W} \boldsymbol{\Sigma}_i\right)^2 +
\sum_{i<l}^k \frac{4}{n_i n_l} \operatorname{tr}\left( \mathbf{W} \boldsymbol{\Sigma}_i \mathbf{W} \boldsymbol{\Sigma}_l\right) \\
& +4 \sum_{i=1}^k \frac{1}{n_i}\left(\sum_{l=1}^k \boldsymbol{\mu}_l-k \boldsymbol{\mu}_i\right)^{\top} \mathbf{W} \boldsymbol{\Sigma}_i \mathbf{W} \left(\sum_{l=1}^k \boldsymbol{\mu}_l-k \boldsymbol{\mu}_i\right).
\end{align*}
Since
\begin{align*}
T_{n}  =&
(k-1) \sum_{i=1}^{k} \frac{1}{n_{i}\left(n_{i}-1\right)} \sum_{j_{1} \neq j_{2}} \left( \mathbf{x}_{i j_{1}} - \boldsymbol{\mu}_{i} + \boldsymbol{\mu}_{i}\right)^{\top}  \mathbf{W} \left( \mathbf{x}_{i j_{2}} - \boldsymbol{\mu}_{i} + \boldsymbol{\mu}_{i}\right)  \\
&-\sum_{i<l}^{k} \frac{2}{n_{i} n_{l}} \sum_{j_{1}, j_{2}} \left( \mathbf{x}_{i j_{1}} - \boldsymbol{\mu}_{i} + \boldsymbol{\mu}_{i}\right)^{\top}  \mathbf{W} \left( \mathbf{x}_{l j_{2}} - \boldsymbol{\mu}_{l} + \boldsymbol{\mu}_{l}\right) \\
=&(k-1) \sum_{i=1}^{k} \frac{1}{n_{i}\left(n_{i}-1\right)} \sum_{j_{1} \neq j_{2}} \left( \mathbf{x}_{i j_{1}} - \boldsymbol{\mu}_{i} \right)^{\top}  \mathbf{W} \left( \mathbf{x}_{i j_{2}} - \boldsymbol{\mu}_{i}\right)  \\
&-\sum_{i<l}^{k} \frac{2}{n_{i} n_{l}} \sum_{j_{1}, j_{2}} \left( \mathbf{x}_{i j_{1}} - \boldsymbol{\mu}_{i}\right)^{\top}  \mathbf{W} \left( \mathbf{x}_{l j_{2}} - \boldsymbol{\mu}_{i} \right)  +  (k-1) \sum_{i=1}^{k} \boldsymbol{\mu}_{i}^{\top} \mathbf{W} \boldsymbol{\mu}_{i} - \sum_{i<l}^{k} 2 \boldsymbol{\mu}_{i}^{\top}  \mathbf{W} \boldsymbol{\mu}_{l} \\ 
&+ (k-1) \sum_{i=1}^{k} \frac{1}{n_{i}\left(n_{i}-1\right)} \sum_{j_{1} \neq j_{2}} \left[  \left( \mathbf{x}_{i j_{1}} - \boldsymbol{\mu}_{i} \right)^{\top}  \mathbf{W} \boldsymbol{\mu}_{i}   + \boldsymbol{\mu}_{i}^{\top} \mathbf{W} \left( \mathbf{x}_{i j_{2}} - \boldsymbol{\mu}_{i} \right)  \right]\\
&- \sum_{i<l}^{k} \frac{2}{n_{i} n_{l}} \sum_{j_{1}, j_{2}} \left[   \left( \mathbf{x}_{i j_{1}} - \boldsymbol{\mu}_{i} \right)^{\top}  \mathbf{W}  \boldsymbol{\mu}_{l} + \boldsymbol{\mu}_{i}^{\top} \mathbf{W} \left( \mathbf{x}_{l j_{2}} - \boldsymbol{\mu}_{l} \right)    \right]  \\
=& \tilde{T}_{n} +  \sum_{i<l}^{k}\left( \boldsymbol{\mu}_{i}- \boldsymbol{\mu}_{l} \right)^{\top} \mathbf{W}\left( \boldsymbol{\mu}_{i}- \boldsymbol{\mu}_{l} \right)  +\hat{T}_{n},
\end{align*}
where  
\begin{align*}
\hat{T}_{n}= & (k-1) \sum_{i=1}^{k} \frac{1}{n_{i}\left(n_{i}-1\right)} \sum_{j_{1} \neq j_{2}} \left[  \left( \mathbf{x}_{i j_{1}} - \boldsymbol{\mu}_{i} \right)^{\top}  \mathbf{W} \boldsymbol{\mu}_{i}   + \boldsymbol{\mu}_{i}^{\top} \mathbf{W} \left( \mathbf{x}_{i j_{2}} - \boldsymbol{\mu}_{i} \right)  \right]\\
&- \sum_{i \neq l}^{k}   \frac{2}{n_{i} } (\sum_{j_{1} } \mathbf{x}_{i j_{1}} - \boldsymbol{\mu}_{i} )^{\top}  \mathbf{W}  \boldsymbol{\mu}_{l}  \\
= & (k-1) \sum_{i=1}^{k} \frac{2}{n_{i}} \sum_{j_{1}=1}^{n_{i}} \left[  \left( \mathbf{x}_{i j_{1}} - \boldsymbol{\mu}_{i} \right)^{\top}  \mathbf{W} \boldsymbol{\mu}_{i}   \right] -     \left[ \sum_{i=1}^{k} \frac{2}{n_{i}}  ( \sum_{j_{1}=1}^{n_{i}} \mathbf{x}_{i j_{1}} - \boldsymbol{\mu}_{i} )^{\top}  \mathbf{W} \right] \left[  \sum_{l=1}^{k}  \boldsymbol{\mu}_{l} \right] \\
&+ \sum_{i=1}^{k} \frac{2}{n_{i}} \sum_{j_{1}=1}^{n_{i}} \left[  \left( \mathbf{x}_{i j_{1}} - \boldsymbol{\mu}_{i} \right)^{\top}  \mathbf{W} \boldsymbol{\mu}_{i}   \right]  \\
= &  k \sum_{i=1}^{k} \frac{2}{n_{i}} \sum_{j_{1}=1}^{n_{i}} \left[  \left( \mathbf{x}_{i j_{1}} - \boldsymbol{\mu}_{i} \right)^{\top}  \mathbf{W} \boldsymbol{\mu}_{i}   \right] - \sum_{i,l=1}^{k} \frac{2}{n_{i}} \sum_{j_{1}=1}^{n_{i}}  \left( \mathbf{x}_{i j_{1}} - \boldsymbol{\mu}_{i} \right)^{\top}  \mathbf{W} \boldsymbol{\mu}_{l}  \\
=& \sum_{i,l=1}^{k} \frac{2}{n_{i}} \sum_{j_{1}=1}^{n_{i}}  \left( \mathbf{x}_{i j_{1}} - \boldsymbol{\mu}_{i} \right)^{\top}  \mathbf{W} \left( \boldsymbol{\mu}_{i} - \boldsymbol{\mu}_{l}\right). 
\end{align*}
Note that
\begin{align*}
\mathbb{E} \left( \mathbf{x}_{i j_{1}} - \boldsymbol{\mu}_{i} \right)^{\top}  \mathbf{W} \left( \boldsymbol{\mu}_{i} - \boldsymbol{\mu}_{l}\right) \left( \boldsymbol{\mu}_{i} - \boldsymbol{\mu}_{l}\right)^{\top} \mathbf{W} \left( \mathbf{x}_{i j_{1}} - \boldsymbol{\mu}_{i} \right) = \operatorname{tr} \left( \boldsymbol{\mu}_{i} - \boldsymbol{\mu}_{l}\right)^{\top} \mathbf{W} \boldsymbol{\Sigma}_{i} \mathbf{W} \left( \boldsymbol{\mu}_{i} - \boldsymbol{\mu}_{l}\right),
\end{align*}
it is straightforward to show that
\begin{align*}
&\operatorname{Var} \left(  \hat{T}_{n} \right)\\
=&  \sum_{i_{1} \neq l_{1}}^{k} \sum_{i_{2} \neq l_{2}}^{k} \frac{4}{n_{i_{1}} n_{i_{2}}} \sum_{j_{1}=1}^{n_{i_{1}}}\sum_{j_{2}=1}^{n_{i_{2}}}  \mathbb{E} \left( \mathbf{x}_{i_{1} j_{1}} - \boldsymbol{\mu}_{i_{1}} \right)^{\top}  \mathbf{W} \left( \boldsymbol{\mu}_{i_{1}} - \boldsymbol{\mu}_{l_{1}}\right)    \left( \boldsymbol{\mu}_{i_{2}} - \boldsymbol{\mu}_{l_{2}}\right)^{\top}  \mathbf{W} \left( \mathbf{x}_{i_{2} j_{2}} - \boldsymbol{\mu}_{i_{2}} \right) \\
=& \sum_{ l_{1} \neq i_{1} = i_{2} \neq l_{2} }^{k} \frac{4}{n_{i_{1}}^{2}  } \sum_{j_{1}=1}^{n_{i_{1}}} \mathbb{E}   \left( \mathbf{x}_{i_{1} j_{1}} - \boldsymbol{\mu}_{i_{1}} \right)^{\top}  \mathbf{W} \left( \boldsymbol{\mu}_{i_{1}} - \boldsymbol{\mu}_{l_{1}}\right)    \left( \boldsymbol{\mu}_{i_{1}} - \boldsymbol{\mu}_{l_{2}}\right)^{\top}  \mathbf{W} \left( \mathbf{x}_{i_{1} j_{1}} - \boldsymbol{\mu}_{i_{1}} \right)                                        \\
&+  \sum_{  i_{1} \neq l_{1} }^{k} \frac{4}{n_{i_{1}}^{2}  } \sum_{j_{1}=1}^{n_{i_{1}}} \mathbb{E}   \left( \mathbf{x}_{i_{1} j_{1}} - \boldsymbol{\mu}_{i_{1}} \right)^{\top}  \mathbf{W} \left( \boldsymbol{\mu}_{i_{1}} - \boldsymbol{\mu}_{l_{1}}\right)    \left( \boldsymbol{\mu}_{i_{1}} - \boldsymbol{\mu}_{l_{1}}\right)^{\top}  \mathbf{W} \left( \mathbf{x}_{i_{1} j_{1}} - \boldsymbol{\mu}_{i_{1}} \right)                                                                           \\
=& 4 \sum_{ l_{1} \neq i_{1} \neq l_{2} }^{k}   \frac{  \left( \boldsymbol{\mu}_{i_{1}} - \boldsymbol{\mu}_{l_{2}}\right)^{\top}  \mathbf{W} \boldsymbol{\Sigma}_{i_{1}}  \mathbf{W} \left( \boldsymbol{\mu}_{i_{1}} - \boldsymbol{\mu}_{l_{1}}\right) }{  n_{i_{1}}^{2}  }   + 4 \sum_{  i_{1} \neq l_{1} }^{k}   \frac{  \left( \boldsymbol{\mu}_{i_{1}} - \boldsymbol{\mu}_{l_{1}}\right)^{\top}  \mathbf{W}  \boldsymbol{\Sigma}_{i_{1}} \mathbf{W} \left( \boldsymbol{\mu}_{i_{1}} - \boldsymbol{\mu}_{l_{1}}\right)         }{  n_{i_{1}}^{2}  } .
\end{align*}
Under Assumption D, we  have
\begin{align*}
\operatorname{Var} \left( \frac{ \hat{T}_{n}  }{ \sigma_{n,k} } \right) = o\left( 1\right).
\end{align*}
\begin{align}
\sigma_{n,k}^{2} =	\sum_{i=1}^k \frac{2(k-1)^2}{n_i\left(n_i-1\right)} \operatorname{tr}\left( \mathbf{W} \boldsymbol{\Sigma}_i\right)^2 +
\sum_{i<l}^k \frac{4}{n_i n_l} \operatorname{tr}\left( \mathbf{W} \boldsymbol{\Sigma}_i \mathbf{W} \boldsymbol{\Sigma}_l\right)
\end{align}	
Then, we can obtain
\begin{align}
\frac{T_{n }-  \sum_{i<l}^{k}\left( \boldsymbol{\mu}_{i}- \boldsymbol{\mu}_{l} \right)^{\top} \mathbf{W}\left( \boldsymbol{\mu}_{i}- \boldsymbol{\mu}_{l} \right)  }{\sigma_{n,k}} =  \frac{ \tilde{T}_{n}   +\hat{T}_{n}   }{\sigma_{n,k}} =  \frac{ \tilde{T}_{n}    }{\sigma_{n,k}}  + o_{p}\left( 1\right).
\end{align}
Therefore, it is sufficient to show that 
\begin{align}
\frac{\tilde{T}_{n}  }{\sigma_{n,k}}	\stackrel{\text { d }}{\longrightarrow} \mathrm{N} (0,1).
\end{align}
To obtain the asymptotic normality of  $\tilde{T}_{n}$, we can assume without loss of generality that $\boldsymbol{\mu}_{1}=\cdots=\boldsymbol{\mu}_{k} =\mathbf{0}$, that is, $\mathbf{x}_{i j}= \boldsymbol{\Gamma}_{i} \mathbf{z}_{i j}$, and 
\begin{align*}
\tilde{T}_{n} =(k-1) \sum_{i=1}^{k} \frac{1}{n_{i}\left(n_{i}-1\right)} \sum_{j_{1} \neq j_{2}}  \mathbf{x}_{i j_{1}}^{\top}  \mathbf{W} \mathbf{x}_{i j_{2}}  -\sum_{i<l}^{k} \frac{2}{n_{i} n_{l}} \sum_{j_{1}, j_{2}}  \mathbf{x}_{i j_{1}} \mathbf{W} \mathbf{x}_{l j_{2}}.
\end{align*}  
Next, we begin to define square integrable martingale. For convenience, let $\mathbf{C}_{j+\sum_{i=1}^{l-1}n_{i}}=\mathbf{W}^{1/2}\mathbf{x}_{lj}$ for $j=1,2,\cdots,n_{l}$.
Define
\begin{align*}
\begin{split}
\phi_{st}= \left \{
\begin{array}{ll}
	\frac{k-1}{n_{l}\left(n_{l}-1\right)}\mathbf{C}_{s}^{\top}\mathbf{C}_{t},                    & s,t\in \Lambda_{l}, l=1,2,\cdots,k\\
	-\frac{1}{n_{l}n_{i}}\mathbf{C}_{s}^{\top}\mathbf{C}_{t},                                 & (s,t)\in \Lambda_{l}\times \Lambda_{i}, 1\leq i <l \leq k
\end{array}
\right.
\end{split}
\end{align*}
where $\Lambda_{l}=\{\sum\limits_{i=1}^{l-1}n_{i}+1,\sum\limits_{i=1}^{l-1}n_{i}+2, \cdots, \sum\limits_{i=1}^{l}n_{i}\}$ with $ l=1, 2, \cdots, k$ and $\sum_{i=1}^{0}n_{i}=0$.

Denote  $V_{nt}=\sum_{s=1}^{t-1}\phi_{st}, S_{nm}=\sum_{t=2}^{m}V_{nt}$,
$\mathcal{F}_{nm}=\sigma\{\mathbf{C}_1,\mathbf{C}_2,\ldots,\mathbf{C}_m\}$ which is the $\sigma$ field  generated by $\{\mathbf{C}_1,\mathbf{C}_2,\ldots,\mathbf{C}_m\}$.
Then, we have
\begin{align*}
\tilde{T}_{n}&=(k-1) \sum_{i=1}^{k} \frac{2}{n_{i}\left(n_{i}-1\right)} \sum_{j_{1}< j_{2}}  \mathbf{x}_{i j_{1}}^{\top}  \mathbf{W} \mathbf{x}_{i j_{2}} -\sum_{i<l}^{k} \frac{2}{n_{i} n_{l}} \sum_{j_{1}, j_{2}}  \mathbf{x}_{i j_{1}}^{\top} \mathbf{W} \mathbf{x}_{l j_{2}}\\
&=2\sum_{l=1}^k\sum_{t=2+\sum_{m=1}^{l-1}n_m}^{\sum_{m=1}^ln_m}\sum_{s=1+\sum_{m=1}^{l-1}n_m}^{t-1}\phi_{st}+2\sum_{l=1}^{k-1}\sum_{s=l+1}^k\sum_{s=1+\sum_{m=1}^{l-1}n_m}^{\sum_{m=1}^ln_m}\sum_{t=1+\sum_{m=1}^{s-1}n_m}^{\sum_{m=1}^sn_m}\phi_{st}\\
&=2\sum_{t=2}^n\sum_{s=1}^{t-1}\phi_{st}=2\sum_{t=2}^n V_{nt}.
\end{align*}
\begin{lemma}\label{matringale}
For each $n$,  $\{S_{nm}, \mathcal{F}_{nm}\}_{m=1}^{n}$is the sequence of zero mean and a square integrable martingale.
\end{lemma}
\textbf{Proof}
Firstly, we can obtain  that $\mathcal{F}_{nt-1}\subseteq\mathcal{F}_{nt}$, for any $1\leq t\leq n$,
and $S_{nm}$ is a square integrable sequence with zero mean.
We only need to show $E(S_{nq}|\mathcal{F}_{nm})=S_{nm}$ for any $q\geq m.$ Note  that if $t \leq m \leq n$,
$$\mathbb{E}( V_{nt}|\mathcal{F} _{nm})=\sum_{s=1}^{t-1}\mathbb{E}(\phi_{st}|\mathcal{F} _{nm})
=\sum_{s=1}^{t-1}\phi_{st}=V_{nt}.$$
If $t> m$,
then $\mathbb{E}( \phi_{st}|\mathcal{F} _{nm}) =c_{st}\cdot \mathbb{E}( \mathbf{C}_s^{\top}\mathbf{C}_t|\mathcal{F} _{nm})$,
where
$c_{st}$ is the coefficient of $\phi_{st}$.  In addition,we will use the symbol $c_{st}$ repeatedly, often without mention.

If $s> m$, as $\mathbf{C}_{s}$ and $\mathbf{C}_{t}$ are both independent of $\mathcal{F}_{nm}$,
$$\mathbb{E}(\phi_{st}|\mathcal{F}_{nm})=\mathbb{E}(\phi_{st})=0.$$
if $s \leq m$
$$\mathbb{E}(\phi_{st}|\mathcal{F}_{n,m})=c_{st} \mathbb{E}(\mathbf{C}_{s}^{\top}\mathbf{C}_{t}|\mathcal{F}_{n,m})
=c_{st} \mathbf{C}_{s}^{\top}\mathbb{E}(\mathbf{C}_{t})=0$$
Hence when $t> m$,
$$\mathbb{E}(V_{nt}|\mathcal{F}_{n,m})=0.$$
In summary, for any $q\geq m$,
\begin{align*}
\mathbb{E}(S_{nq}|\mathcal{F}_{nm})
=\sum_{t=2}^{q}\mathbb{E}(V_{nt}|\mathcal{F}_{nm}) =S_{nm}.
\end{align*}
Then Lemma \ref{matringale} can be established by combining the results above. The proof is complete.

In order to obtain the main results, we employ the martingale central limit theorem (see Corollary 3.1 of \cite{hall1980martingale}). Then Theorem \ref{CLT of main results} can be proved by verifying the convergence of conditional variance and Linderberg condition, which are stated in Lemma \ref{condtional variance} and \ref{the Linderberg condition of maritingale CLT}, respectively.

\begin{lemma}\label{condtional variance}
Under Assumption A-D, as $n,p\to\infty,$  it gets that
\begin{align} \label{condition1 of maritingale CLT}
\frac{\sum_{t=2}^n\mathbb{E}(V_{nt}^2|\mathcal{F}_{nt-1})}{\sigma_{n,k}^{2}}\stackrel{p}\longrightarrow\frac{1}{4}.
\end{align}
\end{lemma}
\textbf{Proof}
For $ \forall t\in \Lambda_{l}$,
\begin{align*}
V_{nt}=\sum_{s=1}^{t-1}\phi_{st} = \sum_{s=1}^{\sum\limits_{i=1}^{l-1}n_{i}+1}\phi_{st} + \sum_{s=\sum\limits_{i=1}^{l-1}n_{i}+1}^{t-1}\phi_{st}= \mathbf{M}_{t-1}^{\top}\mathbf{C}_{t},
\end{align*}
where $\mathbf{M}_{t-1} = \sum_{s=\sum\limits_{i=1}^{l-1}n_{i}+1}^{t-1} \frac{1}{n_{l}\left(n_{l}-1\right)}\mathbf{C}_{s}^{\top} - \sum_{s=1}^{\sum\limits_{i=1}^{l-1}n_{i}+1}\frac{1}{n_{l}n_{i}} \mathbf{C}_{s}^{\top}$.
Thus, we have
\begin{align}
\mathbb{E}(V_{nt}^{2}|\mathcal{F}_{n,t-1})
=&\mathbb{E}\left(\mathbf{M}_{t-1}^{\top}\mathbf{C}_{t}\mathbf{C}_{t}^{\top}\mathbf{M}_{t-1} |\mathcal{F}_{n,t-1}\right)= \mathbf{M}_{t-1}^{\top}\mathbf{W}^{\frac{1}{2}}\tilde{\boldsymbol{\Sigma}}_t \mathbf{W}^{\frac{1}{2}} \mathbf{M}_{t-1}, \label{Expecation of Vnt}
\end{align}
where $\tilde{\boldsymbol{\Sigma}}_t=\boldsymbol{\Sigma}_l$, for $ t\in \Lambda_{l}$. Denote
\begin{align*}
\eta_n=\sum_{t=2}^{n} \mathbb{E}(V_{nt}^2|\mathcal{F}_{n,t-1})  =\sum_{l=1}^{k} \sum_{t\in\Lambda_l} \mathbb{E}(V_{nt}^{2}|\mathcal{F}_{n,t-1}),
\end{align*}
that is, to prove \eqref{condition1 of maritingale CLT}, we just need to demonstrate that
\begin{align}\label{Expectation: condition1 of maritingale CLT}
\mathbb{E}(\eta_n)=\frac{1}{4} \sigma_{n,k}^{2},
\end{align}
and
\begin{align}\label{Variance: condition1 of maritingale CLT}
\operatorname{Var}\left( \eta_n \right) =o(\sigma_{n,k}^{4}).	
\end{align}
by the Chebyshev's inequality. 

According to \eqref{Expecation of Vnt}, we have
\begin{align*}
&\mathbb{E} \left( \sum_{t\in\Lambda_l} \mathbb{E}(V_{nt}^{2}|\mathcal{F}_{n,t-1}) \right) \\
=& \sum_{t\in\Lambda_l} \left( \frac{(k-1)^2}{n_{l}^2 \left(n_{l}-1\right)^2} ( t-1- \sum\limits_{i=1}^{l-1}n_{i} )\operatorname{tr} ( \mathbf{W} \tilde{\boldsymbol{\Sigma}}_t)^{2} +  \sum\limits_{i=1}^{l-1} \frac{1}{n_{i} n_{l}^2} \operatorname{tr}(\mathbf{W} \tilde{\boldsymbol{\Sigma}}_t \mathbf{W} \boldsymbol{\Sigma}_{i} ) \right) \\
=& \frac{(k-1)^2}{2n_{l} \left(n_{l}-1\right)} \operatorname{tr} ( \mathbf{W} \boldsymbol{\Sigma}_{l})^{2} + \sum\limits_{i=1}^{l-1} \frac{1}{n_{i} n_{l}}\operatorname{tr}(\mathbf{W} \boldsymbol{\Sigma}_{l} \mathbf{W} \boldsymbol{\Sigma}_{i} ) .
\end{align*}
Then, it follows that
\begin{align}
\mathbb{E}\left( \eta_n \right) =&\sum_{l=1}^{k} \mathbb{E} \left( \sum_{t\in\Lambda_l} \mathbb{E}(V_{nt}^{2}|\mathcal{F}_{n,t-1}) \right) \notag\\
=& \sum_{l=1}^{k} \frac{(k-1)^2}{2n_{l} \left(n_{l}-1\right)} \operatorname{tr} ( \mathbf{W} \boldsymbol{\Sigma}_{l})^{2} + \sum_{l=1}^{k}\sum\limits_{i=1}^{l-1} \frac{1}{n_{i} n_{l}}\operatorname{tr}(\mathbf{W} \boldsymbol{\Sigma}_{l} \mathbf{W} \boldsymbol{\Sigma}_{i} ) \label{Expectation of eta},
\end{align}
which implies \eqref{Expectation: condition1 of maritingale CLT} is ture.

In what follows,  we consider the variance component of   $\eta_n$. To compute the variance component of $\eta_n$. We can easily
rewrite the term \eqref{Expecation of Vnt} as
\begin{align*}
\mathbb{E}(V_{nt}^{2}|\mathcal{F}_{n,t-1})
=&\mathbb{E}\left(\sum_{s_1,s_2=1}^{t-1}\phi_{s_1 t}\phi_{s_2 t}|\mathcal{F}_{n,t-1}\right)\\
=&\sum_{s_1,s_2=1}^{t-1} c_{s_1t}c_{s_2t}\mathbf{C}_{s_1}^{\top}\mathbb{E}(\mathbf{C}_t\mathbf{C}_t^{\top}|\mathcal{F}_{n,t-1})\mathbf{C}_{s_2}\\
=& \sum_{s_1,s_2=1}^{t-1}c_{s_1t}c_{s_2t}\mathbf{C}_{s_1}^{\top}\mathbb{E}(\mathbf{C}_t\mathbf{C}_t^{\top})\mathbf{C}_{s_2}\\
=& \sum_{s_1,s_2=1}^{t-1}c_{s_1t}c_{s_2t}\mathbf{C}_{s_1}^{\top} \mathbf{W}^{\frac{1}{2}}\tilde{\boldsymbol{\Sigma}}_t \mathbf{W}^{\frac{1}{2}}\mathbf{C}_{s_2}.
\end{align*}
Next, we consider
\begin{align}
\mathbb{E} (\eta_n^{2})=&\mathbb{E}\left[ \sum_{t=2}^{n} \sum_{t\in\Lambda_l}\mathbb{E}(V_{nt}^{2}|\mathcal{F}_{n,t-1})\right]^{2} \notag\\
=&\mathbb{E}\left[ \sum_{t=2}^{n} \sum_{s_1,s_2=1}^{t-1}c_{s_1t}c_{s_2t}
\mathbf{C}_{s_1}^{\top}\mathbf{W}^{\frac{1}{2}}\tilde{\boldsymbol{\Sigma}}_t \mathbf{W}^{\frac{1}{2}}\mathbf{C}_{s_2}\right]^{2} \notag \\
=&	\mathbb{E}(A)+2	\mathbb{E}(B) \label{ moment of two order of eta }
\end{align}

where
\begin{align*}
A=& \sum_{t=2}^{n} \sum_{s_{1},s_{2}=1}^{t-1}\sum_{s_{3},s_{4}=1}^{t-1}c_{s_1t}c_{s_2t}c_{s_3t}c_{s_4t}
\mathbf{C}_{s_{1}}^{\top}\mathbf{W}^{\frac{1}{2}}\tilde{\boldsymbol{\Sigma}}_{t}\mathbf{W}^{\frac{1}{2}}\mathbf{C}_{s_{2}}
\mathbf{C}_{s_{3}}^{\top}\mathbf{W}^{\frac{1}{2}}\tilde{\boldsymbol{\Sigma}}_{t}\mathbf{W}^{\frac{1}{2}}\mathbf{C}_{s_{4}},\\
B=&  \sum_{2 \leq t_1 <  t_2}^{n} \sum_{s_{1},s_{2}=1}^{t_{1}-1}\sum_{s_{3},s_{4}=1}^{t_{2}-1}c_{s_1t_1}c_{s_2t_1}c_{s_3t_2}c_{s_4t_2}
\mathbf{C}_{s_{1}}^{\top}\mathbf{W}^{\frac{1}{2}}\tilde{\boldsymbol{\Sigma}}_{t_{1}}\mathbf{W}^{\frac{1}{2}}\mathbf{C}_{s_{2}}
\mathbf{C}_{s_{3}}^{\top}\mathbf{W}^{\frac{1}{2}}\tilde{\boldsymbol{\Sigma}}_{t_{2}}\mathbf{W}^{\frac{1}{2}}\mathbf{C}_{s_{4}}, 
\end{align*}
Then, we turn to analyze the following term: 
\begin{align*}
\mathbb{E}(A)=  &\sum_{l=1}^{k} \sum_{t\in\Lambda_l}\sum_{s_{1},s_{2}=1}^{t-1}\sum_{s_{3},s_{4}=1}^{t-1}c_{s_1t}c_{s_2t}c_{s_3t}c_{s_4t}
\mathbb{E} \left( \mathbf{C}_{s_{1}}^{\top}\mathbf{W}^{\frac{1}{2}}\tilde{\boldsymbol{\Sigma}}_{t}\mathbf{W}^{\frac{1}{2}}\mathbf{C}_{s_{2}}
\mathbf{C}_{s_{3}}^{\top}\mathbf{W}^{\frac{1}{2}}\tilde{\boldsymbol{\Sigma}}_{t}\mathbf{W}^{\frac{1}{2}}\mathbf{C}_{s_{4}} \right)\\ 
= &\sum_{l=1}^{k} \mathbb{E} A_{l}.
\end{align*}
Without loss of generallity, we need only consider the case $l=3$. We omit it since other terms may be calculated in a similar way.
\begin{align*}
\mathbb{E} A_{3} =&  \sum_{t\in\Lambda_3}\sum_{s_{1},s_{2}=1}^{t-1}\sum_{s_{3},s_{4}=1}^{t-1}c_{s_1t}c_{s_2t}c_{s_3t}c_{s_4t}
\mathbb{E} \left( \mathbf{C}_{s_{1}}^{\top}\mathbf{W}^{\frac{1}{2}}\tilde{\boldsymbol{\Sigma}}_{t}\mathbf{W}^{\frac{1}{2}}\mathbf{C}_{s_{2}}
\mathbf{C}_{s_{3}}^{\top}\mathbf{W}^{\frac{1}{2}}\tilde{\boldsymbol{\Sigma}}_{t}\mathbf{W}^{\frac{1}{2}}\mathbf{C}_{s_{4}} \right) \notag\\
=&O(n^{-8})\sum_{t\in\Lambda_3}\sum_{s_{1},s_{2}=1}^{t-1}\sum_{s_{3},s_{4}=1}^{t-1} \mathbb{E} \left( \mathbf{C}_{s_{1}}^{\top}\mathbf{W}^{\frac{1}{2}}\boldsymbol{\Sigma}_{3}\mathbf{W}^{\frac{1}{2}}\mathbf{C}_{s_{2}}
\mathbf{C}_{s_{3}}^{\top}\mathbf{W}^{\frac{1}{2}}\boldsymbol{\Sigma}_{3}\mathbf{W}^{\frac{1}{2}}\mathbf{C}_{s_{4}} \right).        
\end{align*}
Next, we consider $s_1,s_2,s_3$ and $s_{4}$ in the following three cases: (a) $( s_1= s_2) \neq( s_3= s_4)$ ; (b)$ ( s_1= s_3) \neq( s_2= s_4) $ or $(s_1=s_4)\neq(s_2=s_3)$; (c)$s_1=s_2=s_3=s_4.$

We split the term $\mathbb{E} A_{3}$ into three parts:
\begin{align}\label{the spiliting of the term A}
\mathbb{E} A_{3}=\mathbb{E}(A_{31})+\mathbb{E}(A_{32})+\mathbb{E}(A_{33}),
\end{align}
where
\begin{align}
\mathbb{E}(A_{31}) =&O\left(n^{-8}\right)
\mathbb{E}\left(\sum_{t\in\Lambda_3}\sum_{s_1\neq s_2}^{t-1}
\mathbf{C}_{s_1}^{\top}\mathbf{W}^{\frac{1}{2}}\boldsymbol{\Sigma}_3\mathbf{W}^{\frac{1}{2}}\mathbf{C}_{s_1}
\mathbf{C}_{s_2}^{\top}\mathbf{W}^{\frac{1}{2}}\boldsymbol{\Sigma}_3\mathbf{W}^{\frac{1}{2}}\mathbf{C}_{s_2}\right) \notag \\
=&O(n^{-8})\sum_{t\in\Lambda_3}\sum_{s_1\neq s_2}^{t-1}
\mathbb{E}(\mathbf{C}_{s_1}^{\top}\mathbf{W}^{\frac{1}{2}}\boldsymbol{\Sigma}_3\mathbf{W}^{\frac{1}{2}}\mathbf{C}_{s_1})
\mathbb{E}(\mathbf{C}_{s_2}^{\top}\mathbf{W}^{\frac{1}{2}}\boldsymbol{\Sigma}_3\mathbf{W}^{\frac{1}{2}}\mathbf{C}_{s_2})\notag\\
=&O\left(n^{-5}\right)\sum_{s, t=1}^{3}\operatorname{tr}(\boldsymbol{\Sigma}_3\mathbf{W}\boldsymbol{\Sigma}_{s}\mathbf{W})
\operatorname{tr}(\boldsymbol{\Sigma}_3\mathbf{W}\boldsymbol{\Sigma}_{t}\mathbf{W}),\label{the order of A31}
\end{align}
\begin{align}
\mathbb{E}(A_{32}) =&O\left(n^{-8}\right)
\mathbb{E}\left(\sum_{t\in\Lambda_3}\sum_{s_1\neq s_2}^{t-1}
\mathbf{C}_{s_1}^{\top}\mathbf{W}^{\frac{1}{2}}\boldsymbol{\Sigma}_l\mathbf{W}^{\frac{1}{2}}\mathbf{C}_{s_2}
\mathbf{C}_{s_2}^{\top}\mathbf{W}^{\frac{1}{2}}\boldsymbol{\Sigma}_l\mathbf{W}^{\frac{1}{2}}\mathbf{C}_{s_1}\right)  \notag\\
=&O\left(n^{-5}\right)\sum_{s,t=1}^{3}\operatorname{tr}\mathbf{W} \boldsymbol{\Sigma}_{s}\mathbf{W}\boldsymbol{\Sigma}_3\mathbf{W}\boldsymbol{\Sigma}_{t}\mathbf{W}\boldsymbol{\Sigma}_3,\label{the order of A32}
\end{align}
and	
\begin{align}
\mathbb{E}(A_{33})=&O(n^{-8})
\mathbb{E}\left(\sum_{t\in\Lambda_3}\sum_{s=1}^{t-1}
\mathbf{C}_{s}^{\top}\mathbf{W}^{\frac{1}{2}}\boldsymbol{\Sigma}_3\mathbf{W}^{\frac{1}{2}}\mathbf{C}_{s}
\mathbf{C}_{s}^{\top}\mathbf{W}^{\frac{1}{2}}\boldsymbol{\Sigma}_3\mathbf{W}^{\frac{1}{2}}\mathbf{C}_{s}\right) \notag\\
=&O(n^{-6})\sum_{s=1}^3
\mathbb{E}(\mathbf{z}_{s1}^{\top}\boldsymbol{\Gamma}_{s}^{\top}\mathbf{W}\boldsymbol{\Sigma}_3\mathbf{W}\boldsymbol{\Gamma}_{s}\mathbf{z}_{s1}-\operatorname{tr}\boldsymbol{\Sigma}_3\mathbf{W}\boldsymbol{\Sigma}_{s}\mathbf{W}
+\operatorname{tr}\boldsymbol{\Sigma}_3\mathbf{W}\boldsymbol{\Sigma}_{s}\mathbf{W})^{2} \notag \\
=&O(n^{-6})\sum_{s=1}^3
\left(\mathbb{E}(\mathbf{z}_{s1}^{\top}\boldsymbol{\Gamma}_{s}^{\top}\mathbf{W}\boldsymbol{\Sigma}_3\mathbf{W}\boldsymbol{\Gamma}_{s}\mathbf{z}_{s1}-\operatorname{tr}\boldsymbol{\Sigma}_3\mathbf{W}\boldsymbol{\Sigma}_{s}\mathbf{W})^{2}
+(\operatorname{tr}\boldsymbol{\Sigma}_3\mathbf{W}\boldsymbol{\Sigma}_{s}\mathbf{W})^{2}\right) \notag \\
=&O(n^{-6})\sum_{s=1}^3\left(\operatorname{tr}(\boldsymbol{\Sigma}_3\mathbf{W}\boldsymbol{\Sigma}_s\mathbf{W})^2+(\operatorname{tr}\boldsymbol{\Sigma}_3\mathbf{W}\boldsymbol{\Sigma}_s\mathbf{W})^2\right).\label{the order of A33}
\end{align}
Therefore, by combining \eqref{the spiliting of the term A}, \eqref{the order of A31}, \eqref{the order of A32} and \eqref{the order of A33}, we obtain that
\begin{align*}
\mathbb{E} A_{3}=o(\sigma_{n,k}^{4}).
\end{align*}
Then, we have
\begin{align}
\mathbb{E}(A) = o(\sigma_{n,k}^{4}). \label{the spiliting A of eta}
\end{align}
Using the same method, we can prove that
\begin{align}
2\mathbb{E}(B)
=&\frac{1}{16}\sigma_{n,k}^{4}+o(\sigma_{n,k}^{4}). \label{the spiliting B of eta}
\end{align}
Then, together with \eqref{ moment of two order of eta }, \eqref{the spiliting A of eta} and \eqref{the spiliting B of eta}, we have
\begin{align*}
\mathbb{E}(\eta_n^{2})=\frac{1}{16}\sigma_{n,k}^{4}+o(\sigma_{n,k}^{4}).
\end{align*}
It then follows from \eqref{Expectation of eta} that
\begin{align*}
\operatorname{Var}(\eta_n)=o(\sigma_{n,k}^{4}).
\end{align*}
Consequently,  we have completed the proof \eqref{Variance: condition1 of maritingale CLT}.  Finally, summarizing the above, the proof of \eqref{condition1 of maritingale CLT} is finished.

\begin{lemma}\label{the Linderberg condition of maritingale CLT}
Under Assumptions A-D, as $n,p\to\infty,$  it gets that
\begin{align}\label{condition2 of maritingale CLT}
\sum_{t=2}^{n}\sigma_{n,k}^{-2}\mathbb{E}\{V_{nt}^2I(|V_{nt}|>\epsilon\sigma_{n,k})|\stackrel{p}\longrightarrow 0.
\end{align}
\end{lemma}
\textbf{Proof}	Note that
\begin{align*}
\sum_{t=2}^{n}(\sigma_{n,k})^{-2}\mathbb{E}\{V_{nt}^2I(|V_{nt}|>\epsilon\sigma_{n,k})|\mathcal{F}_{n,t-1}\} 
\leq(\sigma_{n,k})^{-4}\epsilon^{-2}\sum_{t=2}^{n}\mathbb{E}(V_{nt}^4|\mathcal{F}_{n,t-1}).
\end{align*}
Hence, to prove \eqref{condition2 of maritingale CLT}, it is therefore sufficient to demonstrate that
\begin{align} \label{equivalent condition2 of maritingale CLT}
\mathbb{E}\left(\sum_{t=2}^{n}\mathbb{E}(V_{nt}^4|\mathcal{F}_{n,t-1})\right)=o(\sigma_{n,k}^4).
\end{align}
Recall that
$V_{nt}=\sum_{s=1}^{t-1}\phi_{st}, S_{nm}=\sum_{j=2}^{m}V_{nj}$, then
\begin{align*}
\mathbb{E}\left(V_{nt}^4|\mathcal{F}_{n,t-1}\right)= \mathbb{E}\left[(\sum_{s=1}^{t-1}\phi_{st})^4 \mid \mathcal{F}_{n,t-1}\right]
=\mathbb{E}\left[\sum_{s_{1}, s_{2}, s_{3},s_{4} =1}^{t-1} \phi_{s_{1} t} \phi_{s_{2} t} \phi_{s_{3} t} \phi_{s_{4} t}         \mid \mathcal{F}_{n,t-1}\right].
\end{align*}
Additional computations can lead to
\begin{align}
&\mathbb{E}\left(\sum_{t=2}^{n}\mathbb{E}(V_{nt}^4|\mathcal{F}_{n,t-1})\right)\notag \\
=&\sum_{t=2}^{n} \sum_{s_{1}, s_{2}, s_{3},s_{4} =1}^{t-1} \mathbb{E} \left(  \phi_{s_{1} t} \phi_{s_{2} t} \phi_{s_{3} t} \phi_{s_{4} t} \right) \notag \\
=&O(n^{-8}) \sum_{t=2}^{n} \sum_{s_{1}, s_{2}, s_{3},s_{4} =1}^{t-1}
\mathbb{E} \left( \mathbf{C}_t^{\top}\mathbf{C}_{s_{1}} \mathbf{C}_t^{\top}\mathbf{C}_{s_{2}} \mathbf{C}_t^{\top}\mathbf{C}_{s_{3}} \mathbf{C}_t^{\top}\mathbf{C}_{s_{4}}           \right)          \notag     \\
=&O(n^{-8}) \left( 3  \sum_{t=2}^{n}\sum_{s_{1}\neq s_{2}}^{t-1}
\mathbb{E} \left( (\mathbf{C}_t^{\top}\mathbf{C}_{s_{1}} )^2 (\mathbf{C}_t^{\top}\mathbf{C}_{s_{2}})^2 \right)
+\sum_{t=2}^{n}\sum_{s_{1}=1}^{t-1}\mathbb{E}\left( \mathbf{C}_t^{\top}\mathbf{C}_{s_{1}} \right)^4  \right)  \notag \\
=&O(n^{-8}) \left( 3 Q + P\right) , \label{order of 3Q+P}
\end{align}
where
\begin{align*}
Q = \sum_{t=2}^{n}\sum_{s_{1}\neq s_{2}}^{t-1}
\mathbb{E} \left( (\mathbf{C}_t^{\top}\mathbf{C}_{s_{1}} )^2 (\mathbf{C}_t^{\top}\mathbf{C}_{s_{2}})^2 \right)
\end{align*}
and 
\begin{align*}
P= \sum_{t=2}^{n}\sum_{s_{1}=1}^{t-1}\mathbb{E}\left( \mathbf{C}_t^{\top}\mathbf{C}_{s_{1}} \right)^4.
\end{align*}
By applying Proposition A.1. in \cite{chen2010tests}, we have
\begin{align}
Q=&\sum_{t=2}^{n}\sum_{s_{1}\neq s_{2}}^{t-1}
\mathbb{E}(\mathbf{C}_t^{\top}\mathbf{W}^{\frac{1}{2}}\tilde{\boldsymbol{\Sigma}}_{s_{1}} \mathbf{W}^{\frac{1}{2}}\mathbf{C}_t
\mathbf{C}_t^{\top} \mathbf{W}^{\frac{1}{2}}\tilde{\boldsymbol{\Sigma}}_{s_{2}} \mathbf{W}^{\frac{1}{2}}\mathbf{C}_t) \notag \\
=&\sum_{l=1}^{k}\sum_{t\in \Lambda_{l}}\sum_{s_{1}\neq s_{2}}^{t-1}
\mathbb{E}(\mathbf{C}_t^{\top}\mathbf{W}^{\frac{1}{2}}\tilde{\boldsymbol{\Sigma}}_{s_{1}} \mathbf{W}^{\frac{1}{2}}\mathbf{C}_t
\mathbf{C}_t^{\top}\mathbf{W}^{\frac{1}{2}}\tilde{\boldsymbol{\Sigma}}_{s_{2}} \mathbf{W}^{\frac{1}{2}}\mathbf{C}_t) \notag \\
=&o( n^{8}\sigma_{n,k}^4).\label{order of 3Q}
\end{align}
Then, we turn to analyze another term $P$. Note that
\begin{align*}
&\sum_{t=2}^{n_{1} + n_{2}}\sum_{s=1}^{t-1}\mathbb{E}\left( \mathbf{C}_t^{\top}\mathbf{C}_{s} \right)^4 \\
=& \left( \sum_{t=2}^{n_{1}}\sum_{s=1}^{t-1} + \sum_{t=n_{1} +1 }^{n_{1}+ n_{2}}\sum_{s=1}^{n_{1}} + \sum_{t=n_{1} +1 }^{n_{1}+ n_{2}}\sum_{s=n_{1} +1}^{t-1} \right)\mathbb{E}\left( \mathbf{C}_t^{\top}\mathbf{C}_{s} \right)^4\\
=& O(n^{2}) \left( \sum_{l=1}^{2} \left(\operatorname{tr}^2\left(\boldsymbol{\Sigma}_l \mathbf{W}\right)^2
+\operatorname{tr}\left(\boldsymbol{\Sigma}_l \mathbf{W}\right)^4\right) + \operatorname{tr}^2(\boldsymbol{\Sigma}_1 \mathbf{W} \boldsymbol{\Sigma}_2 \mathbf{W}) +\operatorname{tr}(\boldsymbol{\Sigma}_1 \mathbf{W} \boldsymbol{\Sigma}_2 \mathbf{W})^2  \right)
\end{align*}
and 
\begin{align*}
&\sum_{t=n_{1} + n_{2} + 1}^{n_{1} + n_{2} + n_{3}}\sum_{s=1}^{t-1}\mathbb{E}\left( \mathbf{C}_t^{\top}\mathbf{C}_{s} \right)^4 \\
=& \left( \sum_{t=n_{1} + n_{2} + 1}^{n_{1} + n_{2} + n_{3}}\sum_{s=1}^{n_{1} } + \sum_{t=n_{1} + n_{2} + 1 }^{n_{1} + n_{2} + n_{3}}\sum_{s=n_{1} +1}^{n_{1} + n_{2}} + \sum_{t=n_{1} + n_{2} + 1 }^{n_{1} + n_{2} + n_{3}}\sum_{s=n_{1} + n_{2} +1}^{t-1} \right)\mathbb{E}\left( \mathbf{C}_t^{\top}\mathbf{C}_{s} \right)^4\\
=& O(n^{2}) \left(  \operatorname{tr}^2\left(\boldsymbol{\Sigma}_3 \mathbf{W}\right)^2
+\operatorname{tr}\left(\boldsymbol{\Sigma}_3 \mathbf{W}\right)^4 + \sum_{l=1}^{2} \operatorname{tr}^2(\boldsymbol{\Sigma}_3 \mathbf{W} \boldsymbol{\Sigma}_l \mathbf{W}) +\operatorname{tr}(\boldsymbol{\Sigma}_3 \mathbf{W} \boldsymbol{\Sigma}_l \mathbf{W})^2  \right).
\end{align*}
Here, we could decompose the term $P$ as
\begin{align*}
P = \sum_{t=2}^{n}\sum_{s=1}^{t-1}\mathbb{E}\left( \mathbf{C}_t^{\top}\mathbf{C}_{s} \right)^4 
=\sum_{l=1}^{k}\sum_{t\in \Lambda_{l}}\sum_{s=1}^{t-1}\mathbb{E}(\mathbf{C}_s^{\top}\mathbf{C}_t)^4 
=\sum_{l=1}^{k}\sum_{t=\sum\limits_{t=1}^{l-1}n_{t}+1}^{\sum\limits_{t=1}^{l}n_{t}}\sum_{s=1}^{t-1} 
\mathbb{E}(\mathbf{C}_s^{\top}\mathbf{C}_t)^4
\end{align*}
The aforementioned process is then repeated, and in a similar manner, we obtain 
\begin{align}
P 
=&\sum_{l=1}^{k}O(n^2)
\left(\operatorname{tr}^2\left((\boldsymbol{\Sigma}_l  \mathbf{W})^2\right)+\operatorname{tr}\left((\boldsymbol{\Sigma}_l \mathbf{W})^4\right)+\sum_{i=1}^{l-1}\operatorname{tr}^2(\boldsymbol{\Sigma}_i \mathbf{W} \boldsymbol{\Sigma}_l \mathbf{W})
+\sum_{i=1}^{l-1}\operatorname{tr}(\boldsymbol{\Sigma}_i \mathbf{W} \boldsymbol{\Sigma}_l \mathbf{W})^2\right) \notag \\
=&O(n^2)
\left(\sum_{l=1}^{k}\left(\operatorname{tr}^2\left((\boldsymbol{\Sigma}_l \mathbf{W})^2\right)
+\operatorname{tr}\left((\boldsymbol{\Sigma}_l \mathbf{W})^4\right)\right)+\sum_{i<l}^{k} \left( \operatorname{tr}^2(\boldsymbol{\Sigma}_i \mathbf{W}\boldsymbol{\Sigma}_l \mathbf{W})
+\operatorname{tr}(\boldsymbol{\Sigma}_i \mathbf{W}\boldsymbol{\Sigma}_l \mathbf{W})^2 \right)\right). \label{order of P}
\end{align}
By combining \eqref{order of 3Q+P}, \eqref{order of 3Q} and \eqref{order of P}, it follows that
$$\mathbb{E}\left\{\sum_{t=2}^{n}\mathbb{E}(V_{nt}^4|\mathcal{F}_{n,t-1})\right\}=
O(n^{-8})\left( 3Q+P \right)=o(\sigma_{n,k}^4).$$
Then, the proof of \eqref{condition2 of maritingale CLT} is complete.

As a consequence, we have completed the proof of Theorem \ref{CLT of main results}.

\subsection{Proof of Lemma \ref{estimators of trace}}
From the expression of $\hat{\sigma}_{n, k}^{2}$ and $\sigma_{n, k}^{2}$, it suffices to show that
\begin{align*}
\frac{\operatorname{tr} \widehat{ \left(\mathbf{W} \boldsymbol{\Sigma}_{i}\right)^{2} }}{\operatorname{tr}\left(\mathbf{W} \Sigma_{i}\right)^{2}} \stackrel{p}{\longrightarrow} 1,\quad \textnormal{ and }\quad \frac{\operatorname{tr}\left(\widehat{\mathbf{W} \Sigma_{i} \mathbf{W} \Sigma_{l}}\right)}{\operatorname{tr}\left(\mathbf{W} \Sigma_{i} \mathbf{W} \Sigma_{l}\right)} \stackrel{p}{\longrightarrow} 1 .
\end{align*}
The proof of above conclusions is immediate by taking the same procedures in \cite{li2012two}. Thus we omitted the details here.

\bibliographystyle{plainnat-abbrev}
\bibliography{references}  

\begin{thebibliography}{19}
\providecommand{\natexlab}[1]{#1}
\providecommand{\url}[1]{\texttt{#1}}
\expandafter\ifx\csname urlstyle\endcsname\relax
  \providecommand{\doi}[1]{doi: #1}\else
  \providecommand{\doi}{doi: \begingroup \urlstyle{rm}\Url}\fi

\bibitem[Anderson(2003)]{anderson2003introduction}
T.~W. Anderson.
\newblock \emph{An introduction to multivariate statistical analysis}.
\newblock Wiley series in probability and statistics. John Wiley \& Sons,
  Hoboken, NJ, 3 edition, 2003.
\newblock ISBN 0-471-36091-0.

\bibitem[Bai and Saranadasa(1996)]{bai1996effect}
Z.~D. Bai and H.~Saranadasa.
\newblock Effect of high dimension: by an example of a two sample problem.
\newblock \emph{Statist. Sinica}, 6\penalty0 (2):\penalty0 311--329, 1996.
\newblock ISSN 1017-0405.

\bibitem[Cai and Xia(2014)]{cai2014high}
T.~T. Cai and Y.~Xia.
\newblock High-dimensional sparse {MANOVA}.
\newblock \emph{J. Multivariate Anal.}, 131:\penalty0 174--196, 2014.
\newblock ISSN 0047-259X.
\newblock \doi{10.1016/j.jmva.2014.07.002}.
\newblock URL \url{https://doi.org/10.1016/j.jmva.2014.07.002}.

\bibitem[Cai et~al.(2014)Cai, Liu, and Xia]{cai2014two}
T.~T. Cai, W.~D. Liu, and Y.~Xia.
\newblock Two-sample test of high dimensional means under dependence.
\newblock \emph{J. R. Stat. Soc. Ser. B. Stat. Methodol.}, 76\penalty0
  (2):\penalty0 349--372, 2014.
\newblock ISSN 1369-7412.
\newblock \doi{10.1111/rssb.12034}.
\newblock URL \url{https://doi.org/10.1111/rssb.12034}.

\bibitem[Cao et~al.(2019)Cao, Park, and He]{cao2019test}
M.~X. Cao, J.~Park, and D.~J. He.
\newblock A test for the k sample {Behrens}-{Fisher} problem in high
  dimensional data.
\newblock \emph{J. Statist. Plann. Inference}, 201:\penalty0 86--102, 2019.
\newblock ISSN 0378-3758.
\newblock \doi{10.1016/j.jspi.2018.12.002}.
\newblock URL \url{https://doi.org/10.1016/j.jspi.2018.12.002}.

\bibitem[Chen and Qin(2010)]{chen2010two}
S.~X. Chen and Y.~L. Qin.
\newblock A two-sample test for high-dimensional data with applications to
  gene-set testing.
\newblock \emph{Ann. Statist.}, 38\penalty0 (2):\penalty0 808--835, 2010.
\newblock ISSN 0090-5364.
\newblock \doi{10.1214/09-AOS716}.
\newblock URL \url{https://doi.org/10.1214/09-AOS716}.

\bibitem[Chen et~al.(2010)Chen, Zhang, and Zhong]{chen2010tests}
S.~X. Chen, L.~X. Zhang, and P.~S. Zhong.
\newblock Tests for high-dimensional covariance matrices.
\newblock \emph{J. Amer. Statist. Assoc.}, 105\penalty0 (490):\penalty0
  810--819, 2010.
\newblock ISSN 0162-1459.
\newblock \doi{10.1198/jasa.2010.tm09560}.
\newblock URL \url{https://doi.org/10.1198/jasa.2010.tm09560}.

\bibitem[Hall and Heyde(1980)]{hall1980martingale}
P.~Hall and C.~C. Heyde.
\newblock \emph{Martingale limit theory and its application}.
\newblock Probability and mathematical statistics. Academic Press, Inc.
  [Harcourt Brace Jovanovich, Publishers], New York-London, 1980.
\newblock ISBN 0-12-319350-8.

\bibitem[He et~al.(2023)He, Shi, Xu, and Cao]{he2023high}
D.~J. He, H.~J. Shi, K.~Xu, and M.~X. Cao.
\newblock A high-dimensional test for the k-sample {Behrens}-{Fisher} problem.
\newblock \emph{J. Nonparametr. Stat.}, 35\penalty0 (2):\penalty0 239--265,
  2023.
\newblock ISSN 1048-5252.
\newblock \doi{10.1080/10485252.2022.2147172}.
\newblock URL \url{https://doi.org/10.1080/10485252.2022.2147172}.

\bibitem[Himeno and Yamada(2014)]{himeno2014wstimations}
T.~Himeno and T.~Yamada.
\newblock Estimations for some functions of covariance matrix in high dimension
  under non-normality and its applications.
\newblock \emph{J. Multivariate Anal.}, 130:\penalty0 27--44, 2014.
\newblock ISSN 0047-259X.
\newblock \doi{10.1016/j.jmva.2014.04.020}.
\newblock URL \url{https://doi.org/10.1016/j.jmva.2014.04.020}.

\bibitem[Hu et~al.(2017)Hu, Bai, Wang, and Wang]{hu2017on}
J.~Hu, Z.~D. Bai, C.~Wang, and W.~Wang.
\newblock On testing the equality of high dimensional mean vectors with unequal
  covariance matrices.
\newblock \emph{Ann. Inst. Statist. Math.}, 69\penalty0 (2):\penalty0 365--387,
  2017.
\newblock ISSN 0020-3157.
\newblock \doi{10.1007/s10463-015-0543-8}.
\newblock URL \url{https://doi.org/10.1007/s10463-015-0543-8}.

\bibitem[Jiang et~al.(2022)Jiang, Wang, Wen, Jiang, and
  Zhang]{jiang022nonparametric}
Y.~L. Jiang, X.~Q. Wang, C.~H. Wen, Y.~K. Jiang, and H.~P. Zhang.
\newblock Nonparametric two-sample tests of high dimensional mean vectors via
  random integration.
\newblock \emph{J. Amer. Statist. Assoc.}, 0\penalty0 (0):\penalty0 1--14,
  2022.
\newblock \doi{10.1080/01621459.2022.2141636}.
\newblock URL \url{https://doi.org/10.1080/01621459.2022.2141636}.

\bibitem[Li and Chen(2012)]{li2012two}
J.~Li and S.~X. Chen.
\newblock Two sample tests for high-dimensional covariance matrices.
\newblock \emph{Ann. Statist.}, 40\penalty0 (2):\penalty0 908--940, 2012.
\newblock ISSN 0090-5364.
\newblock \doi{10.1214/12-AOS993}.
\newblock URL \url{https://doi.org/10.1214/12-AOS993}.

\bibitem[Schott(2007)]{schott2007some}
J.~R. Schott.
\newblock Some high-dimensional tests for a one-way {MANOVA}.
\newblock \emph{J. Multivariate Anal.}, 98\penalty0 (9):\penalty0 1825--1839,
  2007.
\newblock ISSN 0047-259X.
\newblock \doi{10.1016/j.jmva.2006.11.007}.
\newblock URL \url{https://doi.org/10.1016/j.jmva.2006.11.007}.

\bibitem[Wang et~al.(2015)Wang, Pan, and Li]{wang15}
L.~Wang, B.~Pan, and R.~Z. Li.
\newblock A high-dimensional nonparametric multivariate test for mean vector.
\newblock \emph{J. Amer. Statist. Assoc.}, 110\penalty0 (512):\penalty0
  1658--1669, 2015.
\newblock \doi{10.1080/01621459.2014.988215}.
\newblock URL \url{https://doi.org/10.1080/01621459.2014.988215}.

\bibitem[Xu et~al.(2016)Xu, Lin, Wei, and Pan]{xu2016adaptive}
G.~J. Xu, L.~F. Lin, P.~Wei, and W.~Pan.
\newblock An adaptive two-sample test for high-dimensional means.
\newblock \emph{Biometrika}, 103\penalty0 (3):\penalty0 609--624, 2016.
\newblock ISSN 0006-3444.
\newblock \doi{10.1093/biomet/asw029}.
\newblock URL \url{https://doi.org/10.1093/biomet/asw029}.

\bibitem[Yamada and Himeno(2015)]{yamada2015testing}
T.~Yamada and T.~Himeno.
\newblock Testing homogeneity of mean vectors under heteroscedasticity in
  high-dimension.
\newblock \emph{J. Multivariate Anal.}, 139:\penalty0 7--27, 2015.
\newblock ISSN 0047-259X.
\newblock \doi{10.1016/j.jmva.2015.02.005}.
\newblock URL \url{https://doi.org/10.1016/j.jmva.2015.02.005}.

\bibitem[Zhang et~al.(2020)Zhang, Guo, Zhou, and Cheng]{zhangjt2020simple}
J.~T. Zhang, J.~Guo, B.~Zhou, and M.~Y. Cheng.
\newblock A simple two-sample test in high dimensions based on {L2}-norm.
\newblock \emph{J. Amer. Statist. Assoc.}, 115\penalty0 (530):\penalty0
  1011--1027, 2020.
\newblock ISSN 0162-1459.
\newblock \doi{10.1080/01621459.2019.1604366}.
\newblock URL \url{https://doi.org/10.1080/01621459.2019.1604366}.

\bibitem[Zhou and Guo(2017)]{zhou2017note}
B.~Zhou and J.~Guo.
\newblock A note on the unbiased estimator of $\boldsymbol{\Sigma}^2$.
\newblock \emph{Statist. Probab. Lett.}, 129:\penalty0 141--146, 2017.
\newblock ISSN 0167-7152.
\newblock \doi{10.1016/j.spl.2017.05.014}.
\newblock URL \url{https://doi.org/10.1016/j.spl.2017.05.014}.

\end{thebibliography}
\end{document}